\documentclass[10pt,a4paper]{amsart}
\usepackage{amsmath}
\usepackage{amsfonts}
\usepackage{amsthm}
\usepackage{mathtools}
\usepackage{amscd}
\usepackage[utf8]{inputenc}
\usepackage{tikz-cd}
\usepackage{quiver}
\usepackage{t1enc}
\usepackage[mathscr]{eucal}
\usepackage{indentfirst}
\usepackage{graphicx}
\usepackage{graphics}
\usepackage{pict2e}
\usepackage{enumitem}
\usepackage{epic}
\usepackage{float}
\usepackage{MnSymbol}
\usepackage{multirow}
\usepackage{shuffle} %the symbol shuffle product
\usepackage[table,xcdraw]{xcolor} %background color
\usepackage{hyperref}
\usepackage{tcolorbox}
\numberwithin{equation}{section}
\usepackage[margin=2.9cm]{geometry}
\usepackage{epstopdf} 
\usepackage{cite}

\tikzcdset{scale cd/.style={every label/.append style={scale=#1},
		cells={nodes={scale=#1}}}}

\allowdisplaybreaks

\newmuskip\pFqmuskip
\newcommand*\pFq[6][8]{%
	\begingroup % only local assignments
	\pFqmuskip=#1mu\relax
	\mathchardef\normalcomma=\mathcode`,
	\mathcode`\,=\string"8000
	\begingroup\lccode`\~=`\,
	\lowercase{\endgroup\let~}\pFqcomma
	{}_{#2}F_{#3}{\left[\genfrac..{0pt}{}{#4}{#5};#6\right]}%
	\endgroup
}
\newcommand{\pFqcomma}{{\normalcomma}\mskip\pFqmuskip}

\theoremstyle{plain}
\newtheorem{theorem}{Theorem}[section]
\newtheorem{lemma}[theorem]{Lemma}

\newtheorem{proposition}[theorem]{Proposition}

\theoremstyle{definition}
\newtheorem{definition}[theorem]{Definition}

\newtheorem{remark}[theorem]{Remark}

\newtheorem{example}[theorem]{Example}

\newtheorem{romexample}{Example}

\begin{document}
	
	\title{$q$-analogues of Wilf-Zeilberger seeds and Ramanujan $1/\pi^k$-formulas}
	
	\author[Kam Cheong Au]{Kam Cheong Au}
	
	\address{Rheinische Friedrich-Wilhelms-Universität Bonn \\ Mathematical Institute \\ 53115 Bonn, Germany} 
	
	\email{s6kmauuu@uni-bonn.de}
	\subjclass[2020]{Primary: 11B65, 33D15, 33F10. Secondary: 33C20}
	
	\keywords{Basic hypergeometric series, mock theta function, $q$-analogues, $q$-series, $1/\pi$-series, Wilf-Zeilberger pair}

	\begin{abstract}
	We develop the notion of Wilf-Zeilberger seeds as a powerful framework for generating WZ-pairs and for lifting classical hypergeometric identities to the $q$-setting. As an application, we systematically obtain a large family of $q$-analogues of Ramanujan-type $1/\pi^k$ formulas, extending a body of literature previously limited to sporadic examples. While the majority of these $q$-analogues are modular, we also uncover several curious instances of mock modularity.
	\end{abstract}
	
	\maketitle
		\section{Introduction}

One of Ramanujan’s most impressive discoveries is the celebrated $1/\pi$-series
	\begin{equation}\label{dream_series}\sum_{n\geq 0} \left(\frac{1}{3^8\;11^4}\right)^n \frac{(\frac12)_n (\frac14)_n (\frac34)_n}{(1)_n^3} (26390n+1103) = \frac{9801}{2\sqrt{2}\pi},\end{equation}
here $(a)_n = a(a+1)\cdots (a+n-1)$ is the Pochhammer symbol. Nowadays, the nature of such $1/\pi$-series, namely, \begin{equation}\label{1/pi_gen_series}\sum_{n\geq 0} x^n\frac{\left(\frac12\right)_n \left(\frac{1}{2}-s\right)_n \left(\frac12+s\right)_n}{(1)_n^3} (an+b)  = \frac{\sqrt{c}}{\pi}, \qquad a,b,c,x \in \mathbb{Q}, \quad s\in \left\{0,\frac13,\frac14,\frac16\right\}\end{equation} is considered well-understood: these come from modular parameterization of certain $_3F_2$ series \cite{borwein1987pi, chudnovsky1988approximations, cohen2021rational}. \par
Finding $q$-analogues of these formulas is an active and challenging area of research \cite{guo2018ramanujan, guillera2018wz, wei2020q, guo2020q, campbellhal, guo2019q, hou2019q, chen2021hidden, guo2018q, gasper2011basic}. The difficulty lies in the fact that most such identities cannot be derived from familiar transformation formulas; only a few isolated $q$-instances are currently known, all relying on techniques based on Wilf-Zeilberger pairs. In this article, we substantially extend these techniques and derive a large collection of new $q$-analogues. \par
	Before introducing our methodology, we give the reader an impression of what $q$-analogues of these $1/\pi$-series look like. Throughout, we assume $0<q<1$ and denote the $q$-Pochhammer symbol $(a;q)_n := \prod_{j=0}^{n-1}(1-aq^j)$ and $(a;q)_\infty := \prod_{j=0}^\infty(1-aq^j)$. The $1/\pi$-formula
	\begin{equation}\label{aux_0}\sum_{n\geq 0} \left(\frac{1}{2^6}\right)^n \frac{\left(\frac{1}{2}\right)_n^3}{(1)_n^3} (42n+5)= \frac{16}{\pi}\end{equation}
	 is recovered in the limit $q\to 1$ from each of the following four formulas, which we will derive in Examples \ref{ex_3}, \ref{ex_12}, \ref{ex_18}, \ref{ex_21}:
	\begin{align*}
		&\sum_{n\geq 0} q^{6 n^2}  \frac{\left(3 q^{2 n+1}+3 q^{4 n+2}-2 q^{6 n+1}+q^{6 n+3}-3 q^{8 n+2}-3 q^{10 n+3}+1\right) (q;q^2)_n^6}{(q^{2 n+1}+1)^3 (q^2;q^2)_{2 n}^3} = \frac{(q;q^2)_{\infty }^2}{(q^2;q^2)_{\infty }^2},\\
		&\sum_{n\geq 0} q^{5 n^2} \frac{P_1(n) (q;q^2)_n^2 (q^2;q^4)_n^3}{(q^{4 n+2}+1)^2 (q^4;q^4)_n (q^4;q^4)_{2 n}^2} = \frac{(q^2;q^4)_{\infty }^2}{(q^4;q^4)_{\infty }^2}, \\
		&\sum_{n\geq 0}q^{6 n^2} \frac{ P_2(n)  (q^{6 n+1}+q^{12 n+2}+1) (q^{12};q^{12})_{2 n} (q^3;q^6)_n^4 (q,q^5;q^6)_n^2}{(q^{12 n}+1) (q^{6 n+3}+1) (q^{6 n+1}+1)^2 (q^6;q^6)_{2 n}^2 (q^2,q^{10},q^{12};q^{12})_n^2} =  \frac{(q,q^5,q^7,q^{11};q^{12})_\infty (q^3,q^9;q^{12})_\infty^2}{(q^2,q^6,q^{10},q^{12};q^{12})_\infty^2}, \\
		&\sum_{n\geq 0} \frac{q^{2 n} \left(3+3 q^{1+2 n}-q^{4 n}+2 q^{2+4 n}-3 q^{1+6 n}-3 q^{2+8 n}-q^{3+10 n}\right) (q;q^2)_n^6}{(1+q^{1+2 n})^3 (q^2;q^2)_{2 n}^3} = -\frac{(q;q^2)_{\infty }^4}{q (q^2;q^2)_{\infty }^4} \sum_{n\in \mathbb{Z}} \frac{(-1)^n q^{n^2}}{1-q^{2n-1}},
	\end{align*}
	where {\small $$\begin{aligned}P_1(n) &= 2 q^{4 n+2}+q^{6 n+1}-2 q^{8 n+2}+q^{8 n+4}-q^{10 n+1}+2 q^{10 n+3}-q^{12 n+4}-2 q^{14 n+3}-q^{18 n+5}+1,  \quad \text{ and }\\ P_2(n) &= 1+q^{6 n} (q+q^3)+(q^{12 n}-q^{1+18 n}) (1+q^2+q^4) -q^{2+24 n} (1+q^2)-q^{5+30 n}. \end{aligned}$$}
	
The first three identities are modular, in the sense that their right-hand sides are (up to multiplication by a power of $q$) modular forms in $q=e^{2\pi i \tau}$. The last identity is not modular but so-called mock modular, an extension of classical modularity; to the best of our knowledge, this appears to be the first occurrence of mock modularity in this context. We return to this point later. \par
	
	A \textit{Wilf-Zeilberger (WZ) pair} \cite{AequalsB, wilf1990rational, wilf1992algorithmic} consists of a pair of two-variable functions $F,G$ satisfying
	\begin{equation}\label{WZ}
		F(n+1,k) - F(n,k) = G(n,k+1) - G(n,k).
	\end{equation}
The WZ-method is a powerful tool for discovering and proving hypergeometric series that are otherwise inaccessible by conventional techniques \cite{mohammed2005infinite, pilehrood2010series, guillera2003new, sun2024new}. Its connection to infinite series lies in the fact that a WZ-pair $(F,G)$ yields an identity of the form
	$$\sum_{k\geq 0} F(0,k) = \sum_{n\geq 0} G(n,0) \qquad \text{ (up to simpler terms)}.$$
When $(F,G)$ are products of Pochhammer (resp. $q$-Pochhammer) symbols, one obtains identities between classical (resp. $q$-) hypergeometric series.\par
Iconic results obtained via the WZ-method include the following $1/\pi^2$-series due to Guillera \cite{guillera2003new, guillera2016bilateral}:
		$$\begin{aligned}&\sum_{n\geq 0} \left(\frac{1}{2^4}\right)^n \frac{\left(\frac{1}{2}\right)_n^3 \left(\frac{1}{4}\right)_n \left(\frac{3}{4}\right)_n}{(1)_n^5} (120n^2 + 34n+3) = \frac{32}{\pi^2} \\
		&\sum_{n\geq 0} \left(\frac{-1}{2^{10}}\right)^n \frac{\left(\frac{1}{2}\right)_n^5}{(1)_n^5} (820n^2+180n+13) = \frac{128}{\pi^2}.\end{aligned}$$

These identities appear, like $1/\pi$-formulas, to possess deep arithmetic structure, as suggested by their connections to Calabi-Yau differential equations \cite{almkvist2012ramanujan, zudilin2011arithmetic, zudilin2007ramanujan} and to $L$-functions of Hilbert modular forms \cite{dembele2022special}. However, unlike the case of $1/\pi$-formulas, a concrete general theory remains elusive; to date, the only proofs of these identities rely on the WZ-method. \par

For a given (classical or $q$-) hypergeometric term $F(n,k)$, Gosper’s algorithm decides whether a hypergeometric function $G(n,k)$ exists and constructs it when it does \cite{paule1995mathematica, Koutschan09}. What prevents the WZ-method from becoming a fully developed theory is the lack of a satisfactory and systematic way to construct WZ-pairs. At present, almost all WZ-pairs appearing in the literature are obtained by educated guessing of a candidate $F(n,k)$. \par
	
Inspired by a work of Gessel \cite{gessel1995finding}, the author introduced the notion of a WZ-seed, which allows systematic generation of WZ-pairs. Let $f(a_1,\ldots,a_r)$ be a function on $\mathbb{Z}^r$, and let $\mathbf{S}_{a_j}$ denote the operator that increases the argument $a_j$ by $1$. With this notation, the condition \eqref{WZ} for $(F,G)$ to be a WZ-pair can be written as
$$(\textbf{S}_n-1)F = (\textbf{S}_k-1)G.$$
\begin{definition}\label{WZ_seed_def}
	Let $f(a_1,\cdots,a_m,k)$ be a (classical or $q$-) hypergeometric term in $m+1$ variables. It is called a \textit{WZ-seed} with accessory variables $a_1,\cdots,a_m$ if for each $1\leq j\leq m$, there exists a (classical or $q$-) rational function $r_j$ such that
	$$(\textbf{S}_{a_j}-1)f = (\textbf{S}_k-1)(r_j f).$$
\end{definition}
A WZ-pair is thus simply a WZ-seed with a single accessory variable. The precise connection between WZ-seeds and WZ-pairs is given by the following criterion (Proposition \ref{WZ_crit}).
\begin{proposition}
	The (classical or $q$-) hypergeometric term  $f(a_1,\cdots,a_m,k)$ is a WZ-seed if and only if for all $A_1,\cdots,A_m, K\in \mathbb{Z}$,
	$$F(n,k) = f(a_1+A_1 n,a_2+A_2 n,\cdots,a_m+A_m n,k+Kn)$$
	has a hypergeometric WZ-mate $G(n,k)$, that is, $F(n+1,k) - F(n,k) = G(n,k+1) - G(n,k).$
\end{proposition}

Consequently, it suffices to construct WZ-seeds; a rich source of them arises from well-known hypergeometric summation formulas. We provide an overview of this construction and a (non-exhaustive) list of WZ-seeds in Section 2.2. \par

The power of the WZ-seed concept is illustrated by its use in constructing a previously unknown WZ-pair \cite{au2025wilf}, which establishes the remarkable $1/\pi^3$-formula
\begin{equation}\label{aux_16}\sum_{n\geq 0} \left(\frac{1}{2^6}\right)^n \frac{\left(\frac{1}{2}\right)_n^7}{(1)_n^7} (6 n+1) (28 n^2+8 n+1)  = \frac{32}{\pi^3},\end{equation}
empirically discovered by Gourevich in 2002 and had been conjectural for more than two decades. The second main application, which is the focus of this paper, is that WZ-seeds provide a natural dictionary for translating classical hypergeometric series into $q$-series.

\begin{figure}[h]
	\begin{tikzcd}[column sep=4em, row sep=4em, scale cd=1.2]
	{\text{WZ-seed}} && {q\text{-WZ-seed}} \\
	{\text{WZ-pair}} && {q\text{-WZ-pair}} \\
	{\text{hyp. identity}} && {q\text{-hyp. identity}}
	\arrow["{\text{Section } 2.2}"', curve={height=-30pt}, from=1-1, to=1-3]
	\arrow["\text{Proposition }\ref{WZ_crit}"', from=1-1, to=2-1]
	\arrow["{q\text{ tends to }1}", from=1-3, to=1-1]
	\arrow["\text{Proposition }\ref{WZ_crit}", from=1-3, to=2-3]
	\arrow["\text{Proposition }\ref{WZ_prop}"', from=2-1, to=3-1]
	\arrow["{q\text{ tends to }1}", from=2-3, to=2-1]
	\arrow["\text{Proposition }\ref{WZ_prop}", from=2-3, to=3-3]
	\arrow["{q\text{ tends to }1}", from=3-3, to=3-1]
\end{tikzcd}
\caption{\small A schematic illustration of the three concepts: WZ-seeds, WZ-pair and hyoergeometric summation identity, for both classical and $q$-version. }
\end{figure}

Our flow of ideas is summarized in the diagram. The arrow from the classical to the $q$-world (the curly arrow) exists solely at the level of WZ-seeds (the top row) and is invisible at the level of WZ-pairs or summation identities (the bottom two rows). Our framework is thus not merely a new identity-producing trick, but a structural explanation of why $q$-analogues of many such formulas exist and how to find them. \par

Reaping the fruits of this methodology, we readily obtain a $q$-analogue of the challenging $1/\pi^3$-series \eqref{aux_16} above (Example~\ref{ex_17}):
\begin{multline}\label{1/pi^3}
	\sum_{n\geq 0}\frac{q^{2 n^2} (1-q^{6 n+1}) \left(1+q^{4 n}+q^{2 n+1}-6 q^{6 n+1}+q^{8 n+2}+q^{10 n+1}+q^{12 n+2}\right) (q^4;q^4)_{2 n} (q;q^2)_n^8}{(q^{4 n}+1) (q^{2 n+1}+1) (q^2;q^2)_{2 n}^2 (q^4;q^4)_n^6}
	\\ = \frac{(q;q^2)_{\infty }^4}{(q^2;q^4)_{\infty }^2 (q^4;q^4)_{\infty }^6}.
\end{multline}

Another advantage of the WZ-seed formalism is that it naturally produces identities with additional free parameters. For example, the classical $1/\pi$-formula
$$\sum_{n\geq 0} \left(\frac{1}{9}\right)^n \frac{\left(\frac{1}{2}\right)_n \left(\frac{1}{4}\right)_n \left(\frac{3}{4}\right)_n}{(1)_n^3} (8n+1) = \frac{2\sqrt{3}}{\pi}$$
admits a $q$-analogue
$$\sum_{n\geq 0}\frac{q^{2 n^2} (1-q^{8n+1}) (q;q^2)_{2 n} (q;q^2)_n^2}{(q^2;q^2)_{2 n} (q^6;q^6)_n^2} =\frac{(q;q^2)_{\infty } (q^3;q^6)_{\infty }}{(q^2;q^2)_{\infty } (q^6;q^6)_{\infty }},$$
which in turn admits a two-parameter generalization (Example \ref{ex_16}):
$$\sum_{n\geq0}\frac{q^{2bn+2cn+2n^2}(1-q^{b+4c+8n+1})(q^{b-2c+1},q^{b+4c+1};q^2)_n(q^{-b+2c+1};q^2)_{2n}}{(q^6;q^6)_n(q^{6c+6};q^6)_n(q^{2b+2c+2};q^2)_{2n}}=\frac{(q^{3b+3};q^6)_{\infty}(q^{b+4c+1};q^2)_{\infty}}{(q^{6c+6};q^6){}_\infty(q^{2b+2c+2};q^2)_{\infty}}.$$
In general, the number of free parameters in the resulting identity equals the number of accessory variables (that is, $m$ in Definition~\ref{WZ_seed_def}) of the inducing WZ-seed. \par

As concrete applications of the WZ-seed concept, we construct twenty-one $q$-analogues of $1/\pi^k$-formulas with additional free parameters; all of them appear to be new. By specializing the parameters, we recover all previously known $q$-analogues in the literature and obtain many more. \par

Although WZ-pairs are algebraic in nature, their application to infinite series has some analytic subtleties. We organize our 21 examples according to the analytic complexity involved. In Examples~\ref{ex_1}-\ref{ex_12}, these issues are mild. In Examples~\ref{ex_13}-\ref{ex_18}, we employ Guillera’s technique of flawless WZ-pairs \cite{guillera2025wz} to tame them. These 18 examples yield $q$-analogues of $1/\pi^k$-formulas that are modular\footnote{plainly, right-hand sides of these $q$-analogues are products of $q$-Pochhammer symbols}. For the remaining three examples, the analytic issues are more profound and lead instead to mock modular objects. \par

We already encountered mock modularity in the last $q$-analogue of \eqref{aux_0}, where the series $\sum_{n\in \mathbb{Z}} \frac{(-1)^n q^{n^2}}{1-q^{2n-1}}$ appearing on the right-hand side can be recognized as Zwegers’ $\mu(\tau,\tau;2\tau)$ \cite{zwergers2002}. The remaining two mock modular examples (Examples~\ref{ex_19} and \ref{ex_20}) are
\begin{align*}&\sum_{n\geq 0}\frac{q^{2 n} \left(2-q^{2 n}-q^{4 n+1}\right) (q;q^2)_n^4}{(q^{2 n+1}+1) (q^2;q^2)_{2 n} (q^2;q^2)_n^2} = \frac{(q;q^2)_{\infty }^4}{(q^2;q^2)_{\infty }^3} \omega(q) \qquad \text{ and } \\
&\sum_{n\geq 0}\frac{q^{2 n} \left(3+3 q^{1+2 n}-2 (2-q) q^{1+4 n}-4 q^{1+6 n} (1+q)+2 (1-2 q) q^{1+8 n}+3 q^{2+10 n}+3 q^{3+12 n}\right) (q^2;q^2)_{3 n} (q;q^2)_n^6}{(q^{2 n}-q^n+1) (q^{2 n}+q^n+1) (q^{2 n+1}+1)^3 (q^2;q^2)_n^3 (q^2;q^2)_{2 n}^3} \\ &\quad = \frac{(q;q^2)_{\infty }^6}{(q^2;q^2)_{\infty }^5} \omega(q),
\end{align*}
where $\omega(q) := \sum_{n\geq 1} \frac{q^{n-1}}{(q;q^2)_n^2}$ is Ramanujan’s order~3 mock theta function, essentially Zwegers’ $\mu(3\tau,2\tau;6\tau)$ \cite[Appendix~A]{bringmann2017harmonic}. These two identities are $q$-analogues of
$$\begin{aligned}&\sum_{n\geq 0} \left(\frac{1}{4}\right)^n \frac{\left(\frac{1}{2}\right)_n^3}{(1)_n^3} (6n+1) = \frac{4}{\pi} \\
	&\sum_{n\geq 0} \left(\frac{3^3}{4^3}\right)^n \frac{\left(\frac{1}{2}\right)_n^3 \left(\frac{1}{3}\right)_n \left(\frac{2}{3}\right)_n}{(1)_n^5} (74n^2 + 27n+3) = \frac{48}{\pi^2}.
\end{aligned}$$
All of these mock modular $q$-analogues of $1/\pi^k$-formulas appear to be new. It would therefore be natural to seek further instances. Indeed, our framework predicts the existence of additional examples (see the remark following Example~\ref{ex_21}). Notably, we were unable to find an ordinary modular $q$-analogue of the final $48/\pi^2$ identity above, suggesting that mock modularity may be unavoidable in this case. This observation raises the possibility that $q$-analogues of deeper $1/\pi^k$-series, such as \eqref{dream_series}, may involve higher modular objects. \par

Although we focus chiefly on $q$-analogues of $1/\pi^k$-type formulas, the notion of WZ-seeds is applicable much more broadly. In particular, it can be used to derive $q$-analogues of hypergeometric identities of the form
$$\sum_{n\ge 0} z^n \frac{(a_1)_n\cdots (a_m)_n}{(b_1)_n \cdots (b_m)_n} R(n)
= \text{a simple constant},
\qquad
z\in\mathbb{Q}, \quad 
a_i,b_i\in \mathbb{Q}\cap(0,1], \quad
R(n)\in \mathbb{Q}(n),$$
once a WZ-seed inducing the identity has been identified. As an illustration (Example \ref{zeta3}), we derive a $q$-analogue of the classical identity \cite{amdeberhan1998hypergeometric}
$$-2\zeta(3) =\sum_{n\geq 1}\left(\frac{-1}{2^{10}}\right)^{n} \frac{(1)_n^5}{(\frac12)_n^5}\frac{205 n^2-160 n+32}{n^5},$$
a well-known rapidly convergent series for $\zeta(3)$. Its $q$-analogue takes the form
$$\sum_{k\geq 0} \frac{q^k(1+q^{1+k})}{(1-q^{1+k})^3} = \sum_{n\geq 1} \frac{(-1)^{n-1} q^{5n(n-1)/2} P(n) (q;q)_n^{10}}{(1-q^n)^5 (q;q)_{2n}^5}$$
with $P(n) = 1+5 q^n+10 q^{2 n}+5 q^{-1+3 n} (-1+2 q)+5 (-3+q) q^{-1+4 n}+q^{-2+5 n} (1-24 q+q^2)+(5-15 q) q^{-2+6 n}-5 (-2+q) q^{-2+7 n}+10 q^{-2+8 n}+5 q^{-2+9 n}+q^{-2+10 n}$.

We also note a limitation of our approach. Although it substantially extends the scope of applicability of the WZ method, it does not yield WZ-style proofs, and hence $q$-analogues, for all Ramanujan-type $1/\pi$-formulas. This includes the series \eqref{dream_series}, as well as some comparatively simpler examples such as
$$\begin{aligned}
	\sum_{n\geq 0} \left(\frac{-1}{18^2}\right)^n \frac{\left(\frac12\right)_n \left(\frac14\right)_n \left(\frac34\right)_n}{(1)_n^3} (260n+23) &= \frac{72}{\pi}, \\
	\sum_{n\geq 0} \left(\frac{1}{3^4}\right)^n \frac{\left(\frac12\right)_n \left(\frac14\right)_n \left(\frac34\right)_n}{(1)_n^3} (10n+1) &= \frac{9}{2\sqrt{2}\pi}.
\end{aligned}$$
Moreover, none of the $1/\pi$-formulas in \eqref{1/pi_gen_series} with $s=\frac16$ can currently be derived using WZ techniques. This limitation may indicate the existence of further, as yet undiscovered, WZ-seeds. \par

Although not addressed in the article, a promising application of WZ-seeds would be $p$-adic versions of such $1/\pi^k$-series \cite{ZUDILIN20091848}, most of which are widely open conjectures. For example, the $p$-adic counterpart of \eqref{dream_series} is $$\sum_{n=0}^{p-1} \left(\frac{1}{3^8\;11^4}\right)^n \frac{\left(\frac12\right)_n \left(\frac14\right)_n \left(\frac34\right)_n}{(1)_n^3} (26390n+1103) \stackrel{?}{\equiv} 1103\left(\frac{-2}{p}\right)p \pmod{p^3},$$
whereas for the $1/\pi^3$-formula \eqref{1/pi^3} one expects
$$\sum_{n=0}^{p-1} \left(\frac{1}{2^6}\right)^n \frac{\left(\frac{1}{2}\right)_n^7}{(1)_n^7} (6 n+1) (28 n^2+8 n+1)  \stackrel{?}{\equiv} \left(\frac{-1}{p}\right)p^3 \pmod{p^7}.$$
Both congruences are conjectured to hold for all but finitely many primes\footnote{In these cases, $p>11$.}. The recent technique of creative microscoping by Guo and Zudilin \cite{guo2019q_supercongruence, guo2021dwork, guillera2018wz, guo2020proof, guo2020q, guo2019q}, suggests that such congruences may be accessible by analyzing their $q$-analogues with additional parameters at $p$-th roots of unity. In view of this, the latter congruence may be approachable, while the former remains resistant. \par

~\\[0.02in]
The paper is organized as follows. Section~2 reviews WZ-pairs and WZ-seeds and explains how the translation from classical to $q$-hypergeometric identities occurs at the level of WZ-seeds. Section~3 presents 21 examples of $q$-analogues of $1/\pi^k$-formulas. The first 18 are modular and the last three are mock modular.
	
	%\textbf{Acknowledgements} The author would like to express his sincere gratitude to W. Zudilin, J. Guillera and Z.W. Sun for valuable discussions. A special thank to Wadim for encouraging me to continue mathematical research during an uneasy phase of life. 

	\section{Wilf-Zeilberger seeds}
	\subsection{Definition and basic properties}
	A function $f(n)$ is called a \textit{hypergeometric} (resp. $q$-\textit{hypergeometric}) \textit{term} if $f(n+1)/f(n)$ is a rational function in $n$ (resp. $q^n$), with coefficients in a field $\mathcal{K}$, which is normally clear from the context. Examples of hypergeometric terms are $f(n) = \Gamma(a+n)$ and $(a)_n$, with $\mathcal{K} = \mathbb{Q}(a,n)$; examples of $q$-hypergeometric terms are $f(n) = \Gamma_q(a+n)$ and $(q^a;q)_n$, with $\mathcal{K} = \mathbb{Q}(q,q^a,q^n)$, here $\Gamma_q(z) := \frac{(q;q)_\infty}{(q^z;q)_\infty} (1-q)^{1-z}$ is the $q$-gamma function. A hypergeometric (resp. $q$-hypergeometric) term that is a product of gamma (resp. $q$-gamma) functions is called \textit{proper hypergeometric} \cite{AequalsB}. \par
	A \textit{WZ-pair} $(F(n,k),G(n,k))$ consists of two proper hypergeometric terms satisfying $$F(n+1,k)-F(n,k) = G(n,k+1)-G(n,k).$$
	
	For a proper classical hypergeometric term $F(n,k)$, Gosper's algorithm \cite{gosper1978decision} can find its WZ-mate $G(n,k)$ when it exists, or prove it does not. When such a $G$ exists, $\frac{G(n,k)}{F(n.k)}$ is a rational function in $n,k$, called the \textit{certificate} of the WZ-pair. It is not easy to come up with WZ-pairs: for a random hypergeometric term $F(n,k)$, such a $G$ very likely does not exist. To construct non-trivial WZ-pairs, we require the notion of \textit{WZ-seed}.
	
	Before that, we introduce a few notations. Let $f(a_1,\ldots,a_r)$ be a function $\mathbb{Z}^r \to \mathbb{C}$, and let $\mathbf{S}_{a_j}$ denote the operator that increases the argument $a_j$ by $1$, and also denote $\Delta_{a_j} := \textbf{S}_{a_j}-1.$ The condition for $(F,G)$ to be a WZ-pair can be written as
	$$\Delta_n F = \Delta_k G.$$
	\begin{definition}\label{WZ_seed_def}
		Let $f(a_1,\cdots,a_m,k)$ be a (classical or $q$-) hypergeometric term in $m+1$ variables. It is called a \textit{WZ-seed with accessory variables} $a_1,\cdots,a_m$ if for each $1\leq j\leq m$, there exists a (classical or $q$-) rational function $r_j$ such that
		$\Delta_{a_j} f = \Delta_k(r_j f).$
	\end{definition}
	
	\begin{proposition}\cite[Theorem~3.3]{au2025wilf} \label{WZ_crit}
		A (classical or $q$-) hypergeometric term  $f(a_1,\cdots,a_m,k)$ is a WZ-seed if and only if for all $A_1,\cdots,A_m, K\in \mathbb{Z}$,
		$$F(n,k) = f(a_1+A_1 n,\cdots,a_m+A_m n,k+Kn)$$
		has a hypergeometric WZ-mate.
	\end{proposition}
\typeout{	\begin{proof}
	The implication $\Leftarrow$ is easy: $f(a_1,\cdots,a_m,k)$ is satisfies the latter condition, then there exists a rational function $r_i$ such that
	$$f(a_1,\cdots,a_i+n + 1,\cdots,a_m,k) - f(a_1,\cdots,a_i+n ,\cdots,a_m,k) = \Delta_n (r_i f).$$
	Because the variable $a_i$ always appears with $n$ in the form of $a_i+n$, we have $\Delta_n (r_i f) = \Delta_{a_i} (r_i f)$. Setting $n=0$ gives
	$$\Delta_{a_i}f = f(a_1,\cdots,a_i+ 1,\cdots,a_m,k) - f(a_1,\cdots,a_i ,\cdots,a_m,k) = \Delta_{a_i} (r_i f).$$
	Now we show the other implication. Let $A_i \in \mathbb{Z}, K\in \mathbb{Z}$ be given. For integer $l$, $\textbf{S}_{a_i}^l$ acts by changing $a_i$ to $a_i+l$, the condition means there exists $r_i$ such that
	$$(\textbf{S}_{a_i}-1)f = \Delta_k (r_i f), \qquad 1\leq i\leq m.$$
	In the commutative ring $R = \mathbb{Q}[\textbf{S}_{a_1},\cdots,\textbf{S}_{a_m},\textbf{S}_k,\textbf{S}_{a_1}^{-1},\cdots,\textbf{S}_{a_m}^{-1},\textbf{S}_k^{-1}]$, the ideal $(\textbf{S}_{a_1}-1,\cdots,\textbf{S}_{a_m}-1,\textbf{S}_k-1)$ is maximal, it contains any element of $R$ that vanishes at point $(\textbf{S}_{a_1},\cdots,\textbf{S}_{a_m},\textbf{S}_k)=(1,\cdots,1,1)$. In particular, $\left(\prod_{i=1}^m \textbf{S}_{a_i}^{A_i}\right)\textbf{S}_k^K - 1$ is in this ideal, so there exists $\textbf{g}_i \in R$ and $\textbf{h}\in R$ such that
	$$\left(\prod_{i=1}^m \textbf{S}_{a_i}^{A_i}\right)\textbf{S}_k^K - 1 = \sum_{i=1}^m \textbf{g}_i (\textbf{S}_{a_i}-1) + \textbf{h} (\textbf{S}_k-1),$$
	Applying this to $f(a_1,\cdots,a_m,k)$ gives
	$$f(a_1+A_1,\cdots,a_m+A_m,k+K)-f(a_1,\cdots,a_m,k) = \sum_{i=1}^m \textbf{g}_i \Delta_k(r_i f) + \textbf{h}(\textbf{S}_k-1)f= \Delta_k \left(\sum_{i=1}^m \textbf{g}_i r_i f + \textbf{h}f \right).$$
	Replace $a_i$ by $a_i + A_i n$, $k$ by $k+Kn$, the LHS becomes
	$$\Delta_n (f(a_1+A_1n,\cdots,a_m+A_m n,k+Kn)).$$
	Hence $F(n,k) = f(a_1+A_1n,\cdots,a_m+A_m n,k+Kn)$ has a WZ-mate given by $$(\textbf{g}_i r_i +\textbf{h}) f(a_1+A_1n,\cdots,a_m+A_m n, k+Kn),$$
	completing the proof.
	\end{proof} }

The following empirical observation, due to Gessel in the classical hypergeometric setting \cite{gessel1995finding}, says that a closed-form summation formula induces a WZ-seed.
\begin{tcolorbox}
Let $f(a_1,\cdots,a_m,k)$ be a proper (classical or $q$) hypergeometric term in $a_i$ and $k$. If $$\sum_{k\geq 0} f(a_1,\cdots,a_m,k)$$ exists and is independent of $a_1,\cdots,a_m$, then $f(a_1,\cdots,a_m,k)$ is (very likely) a WZ-seed with accessory parameters $a_1,\cdots,a_m$.
\end{tcolorbox}
Although there are counter-examples to this observation \cite[Example~8.1.1]{AequalsB}, it does hold for almost all classical or $q$-hypergeometric summation formulas found in standard texts \cite{gasper2011basic, slater}.

For example, consider the Gauss $_2F_1$ formula:
$$\pFq{2}{1}{a \quad b}{c}{1} = \sum_{k\geq 0} \frac{(a)_k(b)_k}{(1)_k(c)_k} = \frac{\Gamma(c-a-b)\Gamma(c)}{\Gamma(c-a)\Gamma(c-b)}.$$
Dividing the gamma factor of the RHS into the LHS summand, we arrive at
$$\sum_{k\geq 0} \frac{\Gamma (c-a) \Gamma (a+k) \Gamma (c-b) \Gamma (b+k)}{\Gamma (a) \Gamma (b) \Gamma (k+1) \Gamma (c+k) \Gamma (-a-b+c)} = 1$$
If we denote the summand by $f(a,b,c,k)$, then the above observation would predict it is a WZ-seed with accessory parameters $a,b,c$. Indeed, invoking Gosper's algorithm, one obtains
$$\Delta_a(f) = \Delta_k(r_1 f), \qquad \Delta_b(f) = \Delta_k(r_2 f), \qquad \Delta_c(f) = \Delta_k(r_3 f),$$
with rational functions \begin{equation}\label{aux_12}r_1 = \frac{k (c+k-1)}{a (a-c+1)} ,\qquad 
	r_2 = \frac{k (c+k-1)}{b (b-c+1)}, \qquad
	r_3 = -\frac{k}{a+b-c}.\end{equation}
	
A $q$-analogue of the above summation formula is well-known \cite[(1.5.1)]{gasper2011basic}:
$$\sum_{k\geq 0} \frac{(q^a,q^b;q)_k}{(q^c,q;q)_k} q^{k(c-b-a)} = \frac{(q^{c-a},q^{c-b};q)_\infty}{(q^c,q^{c-a-b};q)_\infty}.$$
Dividing the RHS into the LHS summand again, we arrive at
$$\sum_{k\geq 0} \frac{q^{k (-a-b+c)} \Gamma _q(c-a) \Gamma _q(a+k) \Gamma _q(c-b) \Gamma _q(b+k) }{\Gamma _q(a) \Gamma _q(b) \Gamma _q(k+1) \Gamma _q(c+k) \Gamma _q(-a-b+c)} = 1.$$
Let $f_q(a,b,c,k)$ denote the summand. The above observation again predicts it is a WZ-seed with accessory variables $a,b,c$. Indeed, executing the $q$-Gosper algorithm, we found
$$\Delta_a(f_q) = \Delta_k(r_1 f_q), \qquad \Delta_b(f_q) = \Delta_k(r_2 f_q), \qquad \Delta_c(f_q) = \Delta_k(r_3 f_q),$$
with $q$-rational functions
$$r_1 = \frac{(1-q^k) q^{a-k} (q^{c+k}-q)}{(1-q^a) (q^{a+1}-q^c)}, \quad r_2 = \frac{(1-q^k) q^{b-k} (q^{c+k}-q)}{(1-q^b) (q^{b+1}-q^c)},\quad r_3 = \frac{(1-q^k) q^{a+b}}{q^{a+b}-q^c}.$$
Note that they tend to \eqref{aux_12} when $q\to 1$.

\subsection{List of WZ-seeds}
Here we collect some useful WZ-seeds that we shall employ in our examples. They all have their genesis in summation formulas. The list is by no means exhaustive; seeds that do not generate useful examples are not included. 

We express them in terms of the $q$-gamma function, because when written in this form, their hypergeometric limiting cases under $q\to 1$ become clear, which we denote by removing the $q$ in the subscript:
$$\textsf{Seed}(a,b,\cdots,k) := \lim_{q\to 1}\textsf{Seed}_q(a,b,\cdots,k).$$
These classical versions are also listed in \cite[Section~3.3]{au2025wilf}.

\begin{enumerate}[leftmargin=*]
\item  $\textsf{Gauss2F1}_q(a,b,c,k)$
$$= \frac{q^{k (-a-b+c)} \Gamma _q(c-a) \Gamma _q(a+k) \Gamma _q(c-b) \Gamma _q(b+k) }{\Gamma _q(a) \Gamma _q(b) \Gamma _q(k+1) \Gamma _q(c+k) \Gamma _q(-a-b+c)}$$
Origin: $q$-analogue of $_2F_1$ summation formula \cite[(1.5.1)]{gasper2011basic}

\item  $\textsf{Dixon3F2}_q(a,b,c,k)$
$$ = \frac{ q^{k (a-b-c+1)} (q^{a+k}+1) \Gamma _q(a-b+1) \Gamma _q(a-c+1) \Gamma _q(2 a+k) \Gamma _q(b+k) \Gamma _q(c+k) \Gamma _q(2 a-b-c+1)}{\Gamma _q(a) \Gamma _q(b) \Gamma _q(c) \Gamma _q(k+1) \Gamma _q(a-b-c+1) \Gamma _q(2 a-b+k+1) \Gamma _q(2 a-c+k+1)}$$
Origin: $q$-analogue of Dixon's $_3F_2$ summation formula. \cite[(2.7.2)]{gasper2011basic}

\item $\textsf{Dougall5F4}_q(a,b,c,d,k)$
{\small $$ = \frac{q^{k (a-b-c-d+1)} (1-q^{a+2 k}) \Gamma _q(a+k) \Gamma _q(b+k) \Gamma _q(c+k) \Gamma _q(d+k) \Gamma _q(a-b-c+1) \Gamma _q(a-b-d+1) \Gamma _q(a-c-d+1)}{(1-q) \Gamma _q(b) \Gamma _q(c) \Gamma _q(d) \Gamma _q(k+1)\Gamma _q(a-b+k+1) \Gamma _q(a-c+k+1) \Gamma _q(a-d+k+1) \Gamma _q(a-b-c-d+1)}$$}
Origin: $q$-analogue of very-well-poised $_5F_4$ summation formula \cite[(2.7.1)]{gasper2011basic}

\item $\textsf{Dougall7F6}_q(a,b,c,d,e,k)$
{\small $$ = \frac{\splitfrac{(-1)^{a-e} (1-q^{a+2 k}) q^{\frac{1}{2} (a-e) (a-e+1)+k} \Gamma _q(a+k) \Gamma _q(b+k) \Gamma _q(c+k) \Gamma _q(d+k) \Gamma _q(e+k) \Gamma _q(a-b-c+1)}{\times \Gamma _q(a-b-d+1) \Gamma _q(a-c-d+1) \Gamma _q(-a+b+c+e) \Gamma _q(-a+b+d+e) \Gamma _q(-a+c+d+e) \Gamma _q(2 a-b-c-d-e+k+1)}}{\splitfrac{(1-q) \Gamma _q(b) \Gamma _q(c) \Gamma _q(d) \Gamma _q(e) \Gamma _q(k+1) \Gamma _q(-a+b+e) \Gamma _q(a-b+k+1) \Gamma _q(-a+c+e) \Gamma _q(a-c+k+1)}{\times \Gamma _q(-a+d+e) \Gamma _q(a-d+k+1) \Gamma _q(a-e+k+1) \Gamma _q(a-b-c-d+1) \Gamma _q(2 a-b-c-d-e+1) \Gamma _q(-a+b+c+d+e+k)}}$$}
Origin: $q$-analogue of very-well-poised $2$-balanced terminating $_7F_6$ summation formula \cite[(2.6.2)]{gasper2011basic}

\item $\textsf{Balanced3F2}_q(a,b,c,d,k)$
$$ = \frac{(-1)^{a+b+c-d}q^{\frac{1}{2} (a+b+c-d) (a+b+c-d+1)+k} \Gamma _q(d-a) \Gamma _q(a+k) \Gamma _q(d-b) \Gamma _q(b+k) \Gamma _q(d-c) \Gamma _q(c+k) }{\Gamma _q(a) \Gamma _q(b) \Gamma _q(c) \Gamma _q(k+1) \Gamma _q(d+k) \Gamma _q(-a-b+d) \Gamma _q(-a-c+d) \Gamma _q(-b-c+d) \Gamma _q(a+b+c-d+k+1)}$$ 
Origin: $1$-balanced terminating $_3F_2$ summation formula \cite[(1.7.2)]{gasper2011basic}

\item $\textsf{Watson3F2}_q(a,b,c,k)$
{\small $$=\frac{\splitfrac{q^{k/2-a k-b k+c k} (1+q)^{-2 a-2 b+2 c} (q+q^{1/2+a+b+c+2 k}) \Gamma _q(\frac{1}{2}+a+b-c) \Gamma _{q^2}(\frac{1}{2}-a+c) \Gamma _{q^2}(\frac{1}{2}-b+c) \Gamma _q(2 a+k)}{\times \Gamma _q(2 b+k) \Gamma _{q^2}(\frac{1}{2}+a+b-c+k) \Gamma _{q^2}(c+k) \Gamma _q(\frac{1}{2}+a-b+c+k) \Gamma _q(\frac{1}{2}-a+b+c+k) \Gamma _{q^2}(-\frac{1}{2}+a+b+c+k)}}{ \splitfrac{\Gamma _{q^2}(a) \Gamma _{q^2}(b) \Gamma _{q^2}(\frac{1}{2}+a+b-c) \Gamma _q(\frac{1}{2}-a-b+c) \Gamma _q(1+k) \Gamma _{q^2}(\frac{1}{2}+a+b+k) \Gamma _q(\frac{1}{2}+a+b-c+k)}{\times \Gamma _{q^2}(\frac{1}{2}+a-b+c+k) \Gamma _{q^2}(\frac{1}{2}-a+b+c+k) \Gamma _q(-\frac{1}{2}+a+b+c+k) \Gamma _q(2 c+k)}}$$}
Origin: a $q$-analogue of Watson's $_3F_2$ summation formula \cite[p.~61]{gasper2011basic}

\item $\textsf{Seed1}_q(a,c,k)$
$$ = \frac{q^{-2 a k+2 c k-k} (q+1)^{a-2 c-2 k} \Gamma _q(-a+2 c-1) \Gamma _q(a+2 k)}{\Gamma _q(a) \Gamma _{q^2}(k+1) \Gamma _{q^2}\left(-a+c-\frac{1}{2}\right) \Gamma _{q^2}(c+k)}$$ 
Origin: put $a \to a/2, b\to a/2+1/2$ in $\textsf{Gauss2F1}_q$.

\item $\textsf{Seed2}_q(a,c,d,k)$
$$ = \frac{(-1)^{a+c-d}  q^{a^2+2 a c-2 a d+2 a+c^2-2 c d+2 c+d^2-2 d+2 k} (q+1)^{-2 c-2 k} \Gamma _q(-a+2 d-1) \Gamma _q(a+2 k) \Gamma _{q^2}(d-c) \Gamma _{q^2}(c+k)}{\Gamma _q(a) \Gamma _{q^2}(c) \Gamma _{q^2}(k+1) \Gamma _{q^2}\left(-a+d-\frac{1}{2}\right) \Gamma _{q^2}(d+k) \Gamma _q(-a-2 c+2 d-1) \Gamma _{q^2}\left(a+c-d+k+\frac{3}{2}\right)}$$ 
Origin: put $a \to a/2, b\to a/2+1/2$ in $\textsf{Balanced3F2}_q$.

\item $\textsf{Seed3}_q(a,b,d,k)$
$$ = \frac{(1+q)^{2 d} (1-q^{2 a+4 k})q^{2 a k-2 b k-2 d k+k} \Gamma _{q^2}\left(a-b+\frac{1}{2}\right) \Gamma _{q^2}(a+k) \Gamma _q(b+2 k) \Gamma _{q^2}(d+k)  \Gamma _q(2 a-b-2 d+1)}{(1-q) \Gamma _q(b) \Gamma _{q^2}(d) \Gamma _{q^2}(k+1) \Gamma _{q^2}\left(a-b-d+\frac{1}{2}\right) \Gamma _q(2 a-b+2 k+1) \Gamma _{q^2}(a-d+k+1)}$$ 
Origin: put $b \to b/2, c\to b/2+1/2$ in $\textsf{Dougall5F4}_q$. 

\item $\textsf{Seed4}_q(a,b,k)$
$$ = \frac{(q^2+q+1)^{b+1} (1-q^{3 a+6 k}) q^{3 a k-3 b k} \Gamma _q(3 a-2 b) \Gamma _{q^3}(a+k) \Gamma _q(b+3 k)}{(1-q^3) \Gamma _q(b) \Gamma _{q^3}(k+1) \Gamma _{q^3}(a-b) \Gamma _q(3 a-b+3 k+1)}$$ 
Origin: put $b \to b/3, c\to b/3+1/3, d\to b/3+2/3$ in $\textsf{Dougall5F4}_q$. 

\item $\textsf{Seed7}_q(a,b,d,e,k)$
{\small $$ = \frac{\splitfrac{(-1)^{a-e} q^{a^2-2 a e+a+e^2-e+2 k} (1-q^{2 a+4 k}) \Gamma _{q^2}\left(a-b+\frac{1}{2}\right) \Gamma _{q^2}(a+k) \Gamma _q(b+2 k)}{\times \Gamma _{q^2}(d+k) \Gamma _{q^2}(e+k) \Gamma _q(2 a-b-2 d+1) \Gamma _{q^2}\left(-a+b+e+\frac{1}{2}\right) \Gamma _q(-2 a+b+2 d+2 e) \Gamma _{q^2}\left(2 a-b-d-e+k+\frac{1}{2}\right)}}{\splitfrac{(1-q^2) \Gamma _q(b) \Gamma _{q^2}(d) \Gamma _{q^2}(e) \Gamma _{q^2}(k+1) \Gamma _{q^2}\left(a-b-d+\frac{1}{2}\right) \Gamma _q(-2 a+b+2 e)\Gamma _q(2 a-b+2 k+1) \Gamma _{q^2}(-a+d+e)}{\times \Gamma _{q^2}(a-d+k+1) \Gamma _{q^2}(a-e+k+1) \Gamma _{q^2}\left(2 a-b-d-e+\frac{1}{2}\right) \Gamma _{q^2}\left(-a+b+d+e+k+\frac{1}{2}\right)}} $$}
Origin: put $b \to b/2, c\to b/2+1/2$ in $\textsf{Dougall7F6}_q$.

\item $\textsf{Seed8}_q(a,b,e,k)$
{\small $$ = \frac{\splitfrac{(-1)^{a-e} q^{\frac{3 a}{2}+\frac{3 a^2}{2}-\frac{3 e}{2}-3 a e+\frac{3 e^2}{2}+3 k} (1-q^{3 a+6 k}) \Gamma _q(3 a-2 b) \Gamma _q(1-3 a+2 b+3 e)}{\times \Gamma _{q^3}(a+k) \Gamma _{q^3}(2 a-b-e+k) \Gamma _{q^3}(e+k) \Gamma _q(b+3 k)}}{ \splitfrac{(1-q^3) \Gamma _{q^3}(a-b) \Gamma _q(b) \Gamma _{q^3}(2 a-b-e) \Gamma _{q^3}(e) \Gamma _q(-3 a+b+3 e)}{\times \Gamma _{q^3}(1+k) \Gamma _{q^3}(1+a-e+k) \Gamma _{q^3}(1-a+b+e+k) \Gamma _q(1+3 a-b+3 k)}}$$}
Origin: put $b \to b/3, c\to b/3+1/2, d\to b/3+2/3$ in $\textsf{Dougall7F6}_q$. 

\item $\textsf{Seed9}_q(a,b,d,k)$
$$ = \frac{\splitfrac{(-1)^b  q^{2 a b-b^2-b d+2 k}  \left(1-q^{2 (a+2 k)}\right) (q+1)^{4 a-2 b-2 d} \Gamma _{q^2}\left(a-b+\frac{1}{2}\right) \Gamma _{q^2}(a+k) \Gamma _q(b+2 k)}{ \Gamma _q(d+2 k)\Gamma _q(2 a-b-d+1) \Gamma _q(-2 a+2 b+d+1) \Gamma _q(-2 a+b+2 d+1) \Gamma _{q^2}(2 a-b-d+k)}}{\splitfrac{(1-q) \Gamma _q(b) \Gamma _q(d) \Gamma _{q^2}(k+1) \Gamma _{q^2}\left(-a+d+\frac{1}{2}\right) \Gamma _{q^2}(2 a-b-d)}{\Gamma _q(-2 a+b+d+1) \Gamma _q(2 a-b+2 k+1) \Gamma _q(2 a-d+2 k+1) \Gamma _{q^2}(-a+b+d+k+1)}}$$ 
Origin: put $b \to b/2, c\to b/2+1/2, d\to d/2, e\to d/2+1/2$ in $\textsf{Dougall7F6}_q$. 

\item $\textsf{Seed10}_q(a,d,k)$
$$ = \frac{(-1)^{a-d} (q^2+q+1)^{-a-3 k} q^{\frac{3 a^2}{2}-3 a d+\frac{9 a}{2}+\frac{3 d^2}{2}-\frac{9 d}{2}+3 k} \Gamma _q(-a+3 d-2) \Gamma _q(a+3 k)}{\Gamma _q(a) \Gamma _{q^3}(k+1) \Gamma _q(-2 a+3 d-3) \Gamma _{q^3}(d+k) \Gamma _{q^3}(a-d+k+2)}$$ 
Origin: put $a\to a/3, b\to a/3+1/3, c\to a/3+2/3$ in $\textsf{Balanced3F2}_q$. 
\end{enumerate}

All entries in the list are WZ-seeds. To verify this, one requires an implementation of Gosper’s algorithm. Several such implementations exist in various programming languages. The ones used by the author were developed at the Research Institute for Symbolic Computation (RISC) and are written in the \textit{Mathematica} programming language:

\begin{itemize}[leftmargin=*]
	\item \texttt{fastZeil}, by P.~Paule, M.~Schorn, and A.~Riese \cite{paule1995mathematica}, for the classical Gosper algorithm;
	\item \texttt{qZeil}, by A.~Riese \cite{paule1995mathematica}, for the $q$-Gosper algorithm;
	\item \texttt{HolonomicFunctions}, by C.~Koutschan \cite{Koutschan09}, which includes Gosper’s algorithm as a special case of creative telescoping for both classical and $q$-hypergeometric terms. While more general in scope, this package is typically slower than the two specialized packages above.
\end{itemize}

With the help of an implementation of Gosper's algorithm in these packages, one computes the $q$-rational function $r_i$ in the definition of a WZ-seed $$\Delta_{a_i} f = \Delta_k(r_i f).$$ Once these $r_i$ have been computed, the proof of Proposition \ref{WZ_crit} (see \cite{au2025wilf}) then provides an efficient way to construct the $G(n,k)$ for $$F(n,k) = f(a_1+A_1 n,\cdots,a_m+A_m n,k+Kn), \qquad A_1,\cdots,A_m,K\in \mathbb{Z},$$ without further use of Gosper’s algorithm\footnote{In principle, one could still apply Gosper’s algorithm to determine the WZ-mate of $F(n,k)$ directly, but this quickly becomes computationally infeasible even for small values of $A_1,\ldots,A_m,K$}. See the end of this section for an implementation of this procedure. Rather than invoking Gosper’s algorithm, it relies solely on the precomputed functions $r_i$ and suffices to determine all WZ-mates required in this article.

\subsection{Two introductory examples}
The following proposition links WZ-pairs and infinite series. 
\begin{proposition}\cite{au2025wilf}\label{WZ_prop}
	Suppose $F,G: \mathbb{N}^2\to \mathbb{C}$ are two functions such that
	$$F(n+1,k) - F(n,k) = G(n,k+1)-G(n,k)$$
	and the following conditions hold:\begin{itemize}[leftmargin=*]
		\item $\sum_{k\geq 0} F(0,k)$ converges;
		\item $\sum_{n\geq 0} G(n,0)$ converges;
		\item $\lim_{k\to\infty} G(n,k) = g(n)$ exists for each $n\in \mathbb{N}$ and $\sum_{n\geq 0} g(n)$ converges;
	\end{itemize}
	then $\lim_{n\to \infty} \sum_{k\geq 0} F(n,k)$ exists and is finite, also
	$$\sum_{k\geq 0} F(0,k) + \sum_{n\geq 0} g(n) = \sum_{n\geq 0} G(n,0) + \lim_{n\to \infty} \sum_{k\geq 0} F(n,k).$$
\end{proposition}
\typeout{\begin{proof}[Proof sketch]
Applying $\sum_{k=0}^{K-1}\sum_{n=0}^{N-1}$ on both sides of $F(n+1,k) - F(n,k) = G(n,k+1)-G(n,k)$, then let $K\to \infty$ and then $N\to \infty$.
\end{proof}}

Throughout this article, we tacitly assume $0<q<1$ whenever $q$-series are considered. The WZ-seeds given above are expressed using the $q$-gamma function, which has poor analytical properties. In practice, before passing to infinite series, $(F,G)$ must first be converted to a certain "normal form". More precisely, this means (here $C$ is independent of $n$ and $k$)
\begin{itemize}[leftmargin=*]
\item for a classical hypergeometric term, a normal form is 
$$(\text{a rational function in }n,k) \times \frac{\text{product of Pochhammer symbols}}{\text{product of Pochhammer symbols}},$$
where each Pochhammer symbol occurring is of the form  \begin{itemize}
	\item $(C+Kk)_{Nn}, K\in \mathbb{Z}, N\in \mathbb{N}$ or
	\item $(C)_{Kn}, K\in \mathbb{N}.$
	\end{itemize}
\item for a $q$-hypergeometric term, a normal form is 
$$(\text{a $q$-rational function in }n,k) \times \frac{\text{product of $q$-Pochhammer symbols}}{\text{product of $q$-Pochhammer symbols}},$$
where each $q$-Pochhammer symbol occurring is of the form  \begin{itemize}
	\item $(q^{C+Kk};q^r)_{Nn}, K\in \mathbb{Z}, N\in \mathbb{N}, r\in \mathbb{N}$ or
	\item $(q^{C};q^r)_{Nk}, K\in \mathbb{N}, r\in \mathbb{N}.$
\end{itemize}
\end{itemize}
In short, subscripts of Pochhammer symbols in normalized form must have a positive coefficient in $n$ or $k$; a term like $(q^{a+k};q)_{-n}$ is not allowed. In the definition of WZ-pair
$$F(n+1,k) - F(n,k) = G(n,k+1) - G(n,k),$$
multiplying $(F,G)$ by a 1-periodic function in $n,k$ preserves the WZ-property. With the help of the following lemma, one can convert any product of gamma (or $q$-gamma) functions to normal form.
\begin{lemma}
(a) The function $\Gamma(x)\Gamma(1-x) (-1)^x$ is $1$-periodic in $x$. \\
(b) The function $\Gamma_q(x)\Gamma_q(1-x) (-1)^x q^{-x^2/2+x/2}$ is $1$-periodic in $x$. 
\end{lemma}
\begin{proof}
We give a "purely algebraic" proof. For the first one, write $f(x) = (-1)^x \Gamma(x)\Gamma(1-x)$. Using $\Gamma(1+x) = x \Gamma(x)$, one sees
$$\frac{f(x+1)}{f(x)} = -\frac{\Gamma(1+x)}{\Gamma(x)} \frac{\Gamma(-x)}{\Gamma(1-x)} = -x \frac{1}{-x} = 1,$$
so $f(x)$ is $1$-periodic. \par
For the second one, write $f(x) = \Gamma_q(x)\Gamma_q(1-x) (-1)^x q^{-x^2/2+x/2}$, using $\Gamma_q(1+x) = \frac{1-q^x}{1-q} \Gamma_q(x)$, a similar computation shows $f(x+1)/f(x) = 1$, so $f(x)$ is $1$-periodic.
\end{proof}

For a $q$-gamma factor $\Gamma_q(x)$, if its coefficient of $n$ in $x = C+Kk+Nn$ is positive, i.e. $N>0$, then one can directly convert it to the $q$-Pochhammer symbol via $$\Gamma_q(x) = \Gamma_q(C+Kk+Nn) = (q^{C+Kk};q)_{Nn} (1-q)^{Nn} \Gamma_q(C + Kk);$$
whereas if $N<0$, then first replace $\Gamma_q(x)$ by $$\frac{q^{x^2/2-x/2} (-1)^x}{\Gamma_q(1-x)}$$ and then repeat the above replacement process. In any case, one sees quite easily that a normal form can be obtained.

Note that the certificate of $(F,G)$, which is the rational function $R(n,k)=\frac{G(n,k)}{F(n,k)}$, is invariant under normalization, since both $F$ and $G$ are multiplied by the same $1$-periodic functions. 

\begin{example}\label{zeta3}
	Start from $$\textsf{Dougall5F4}(-n,-n,-n,-n,k+2 n+1) = \frac{(2 (k+2 n+1)-n) \Gamma (n+1)^3 \Gamma (k+n+1)^4}{\Gamma (-n)^3 \Gamma (2 n+1) \Gamma (k+2 n+2)^4}.$$
	Replacing $\Gamma(-n)$ by $\frac{(-1)^n}{\Gamma(1+n)}$, we normalize the gamma product into
	$$F(n,k) := \frac{(-1)^n (2 (k+2 n+1)-n) (1)_n^6 (k+1)_n^4}{(k+2 n+1)^4 (1)_{2 n} (k+1)_{2 n}^4}.$$
	Then, for the certificate $R(n,k)$ of this WZ-pair, $2 (2 n+1) (k+2 n+2)^4 (2 k+3 n+2)R(n,k)$ equals {\small $2 k^6+k^5 (26 n+22)+k^4 (141 n^2+240 n+101)+8 k^3 (51 n^3+131 n^2+111 n+31)+8 k^2 (n+1)^2 (83 n^2+120 n+43)+2 k (n+1)^3 (287 n^2+382 n+127)+(n+1)^4 (205 n^2+250 n+77)$}. Writing $G(n,k) = F(n,k) R(n,k)$ and applying Proposition \ref{WZ_prop}, we find  that\footnote{since this is only an illustrative example on classical hypergeometric case, we do not elaborate on these points, see \cite[Section~3.5]{au2025wilf} for the rigorous justification.} $g(n) = 0$ and $\lim_{n\to \infty} \sum_{k\geq 0}F(n,k) = 0$, so we have $$\sum_{k\geq 0} F(0,k) = \sum_{n\geq 0} G(n,0) \iff \sum_{k\geq 0}\frac{2}{(1+k)^3} = \sum_{n\geq 0} \frac{(-1)^n (205 n^2+250 n+77) (1)_n^{10}}{32 (2 n+1)^5 (1)_{2 n}^5},$$ giving a famous fast-converging formula of $\zeta(3)$ due to \cite{amdeberhan1998hypergeometric}:
	\begin{equation}\label{aux_13}-2\zeta(3) =\sum_{n\geq 1}\left(\frac{-1}{2^{10}}\right)^{n} \frac{(1)_n^5}{(\frac12)_n^5}\frac{205 n^2-160 n+32}{n^5}.\end{equation}
	How to derive a $q$-analogue of this formula? The answer is simply repeat the above process but with the $q$-version of the corresponding WZ-seed, that is
	$$\textsf{Dougall5F4}_q(-n,-n,-n,-n,k+2 n+1) = \frac{q^{2 k n+k+4 n^2+4 n+1} (1-q^{2 k+3 n+2}) \Gamma _q(n+1)^3 \Gamma _q(k+n+1)^4}{(1-q) \Gamma _q(-n)^3 \Gamma _q(2 n+1) \Gamma _q(k+2 n+2)^4}.$$
	Replacing $\Gamma_q(-n)$ by $\frac{(-1)^n q^{n(n-1)/2}}{\Gamma_q(1+n)}$, we normalize the gamma product into $$F(n,k) := \frac{(-1)^n q^{2 k n+k+\frac{5 n^2}{2}+\frac{5 n}{2}+1} (1-q^{2 k+3 n+2}) (q;q)_n^6 (q^{k+1};q)_n^4}{(1-q^{k+2 n+1})^4 (q;q)_{2 n} (q^{k+1};q)_{2 n}^4},$$
	Then for the certificate $R(n,k)$ of this WZ-pair, $$(q^{n+1}+1) (q^{2 n+1}-1) (q^{k+2 n+2}-1)^4 (q^{2 k+3 n+2}-1) R(n,k)$$ equals {\small $1+q^{1+n}-4 q^{2+k+2 n}-4 q^{10+5 k+12 n}+q^{11+6 k+13 n}+q^{12+6 k+14 n}+q^{8+4 k+10 n} (6+5 q)+q^{3+2 k+4 n} (5+6 q)+q^{7+4 k+9 n} (6+6 q+q^2)+q^{3+2 k+5 n} (1+6 q+6 q^2)-q^{2+k+3 n} (4+4 q+q^k)-q^{9+4 k+11 n} (q+4 q^k+4 q^{1+k})+q^{4+2 k+6 n} (1+10 q-4 q^{2+k}+q^{1+2 k})-q^{5+2 k+7 n} (5 q+4 q^k+24 q^{1+k}+4 q^{2+k}+5 q^{1+2 k})+q^{6+2 k+8 n} (q-4 q^k+10 q^{1+2 k}+q^{2+2 k})$}. 
	Applying Proposition \ref{WZ_prop}, it is easy to see (we assume, as always, that $0<q<1$)
	$$g(n) = \lim_{k\to \infty} G(n,k) = 0, \qquad \lim_{n\to \infty} \sum_{k\geq 0}F(n,k) = 0,$$
	because of the exponentially decaying factor $q^{2kn+k+5n^2/2}$ in the numerator\footnote{here the $q$-case is easier than the classical case: the vanishing of $\lim_{n\to \infty} \sum_{k\geq 0}F(n,k)$ is immediate for the $q$-case but requires some asymptotic analysis for the classical case}. The equality $\sum_{k\geq 0} F(0,k) = \sum_{n\geq 0} G(n,0)$ yields
	$$\sum_{k\geq 0} \frac{q^k(1+q^{1+k})}{(1-q^{1+k})^3} = \sum_{n\geq 1} \frac{(-1)^{n-1} q^{5n(n-1)/2} P(n) (q;q)_n^{10}}{(1-q^n)^5 (q;q)_{2n}^5}$$
	with $P(n) = 1+5 q^n+10 q^{2 n}+5 q^{-1+3 n} (-1+2 q)+5 (-3+q) q^{-1+4 n}+q^{-2+5 n} (1-24 q+q^2)+(5-15 q) q^{-2+6 n}-5 (-2+q) q^{-2+7 n}+10 q^{-2+8 n}+5 q^{-2+9 n}+q^{-2+10 n}$. This is a $q$-analogue of equation \eqref{aux_13}.
\end{example}

Next we give an example with accessory parameters.
\begin{example}\label{gauss_intro_ex}
	Let $a,b,c,d$ be in a neighbourhood of $0$. We start with a hypergeometric WZ-seed $$\textsf{Gauss2F1}(a-d-n,b-d-n,c-d+1,d+k+2 n+1).$$
	Normalizing, write $$F(n,k) :=\frac{(a+1)_k (b+1)_k (-a+c+1)_n (-a+d+1)_n (a+k+1)_n (-b+c+1)_n (-b+d+1)_n (b+k+1)_n}{(c+k+2 n+1)  (d+k+2 n+1)  (c+1)_k(d+1)_k (c+k+1)_{2 n} (d+k+1)_{2 n} (-a-b+c+d+1)_{2 n}}.$$
	One computes the certificate $R(n,k) = \frac{G(n,k)}{F(n,k)}$ (which we do not display here). Applying Proposition \ref{WZ_prop}, it can be shown that\footnote{since this is only an illustrative example in the classical hypergeometric case, we do not elaborate on these points, see \cite[Section~3.5]{au2025wilf} for rigorous justification.} $g(n) = 0$ and $\lim_{n\to \infty} \sum_{k\geq 0}F(n,k) = 0$. The equality $\sum_{k\geq 0} F(0,k) = \sum_{n\geq 0} G(n,0)$ is explicitly
	\begin{multline} \label{ex1} \sum _{k\geq 0} \frac{(1+a)_k (1+b)_k}{(1+c)_{k+1} (1+d)_{k+1}} \\  = \sum_{n\geq 1} \frac{ (1+a)_{n-1} (1+b)_{n-1} (1-a+c)_{n-1} (1-b+c)_{n-1} (1-a+d)_{n-1} (1-b+d)_{n-1} P(n)}{(1+c)_{2 n} (1+d)_{2 n} (1-a-b+c+d)_{2 n}}\end{multline}
	here {\small $P(n) = a^2 b^2-a^2 b c-a b^2 c+a b c^2-a^2 b d-a b^2 d+a c d+a^2 c d+b c d+3 a b c d+b^2 c d-c^2 d-2 a c^2 d-2 b c^2 d+c^3 d+a b d^2-c d^2-2 a c d^2-2 b c d^2+2 c^2 d^2+c d^3-2 a^2 b n-2 a b^2 n+2 a c n+a^2 c n+2 b c n+6 a b c n+b^2 c n-2 c^2 n-3 a c^2 n-3 b c^2 n+2 c^3 n+2 a d n+a^2 d n+2 b d n+6 a b d n+b^2 d n-6 c d n-11 a c d n-11 b c d n+10 c^2 d n-2 d^2 n-3 a d^2 n-3 b d^2 n+10 c d^2 n+2 d^3 n+4 a n^2+a^2 n^2+4 b n^2+8 a b n^2+b^2 n^2-8 c n^2-13 a c n^2-13 b c n^2+13 c^2 n^2-8 d n^2-13 a d n^2-13 b d n^2+29 c d n^2+13 d^2 n^2-8 n^3-14 a n^3-14 b n^3+28 c n^3+28 d n^3+21 n^4$}. When $a=b=c=d=0$, we recover an identity due to Zeilberger \cite{wilf1992algorithmic}:
	\begin{equation} \label{ex1S}\zeta(2) = \sum_{n\geq 1} \left(\frac{1}{2^6}\right)^n \frac{(1)_n^3}{(\frac12)_n^3} \frac{21n-8}{n^3}. \end{equation}
	
	To derive a $q$-version of \eqref{ex1}, one simply uses the corresponding $q$-version of the WZ-seed: $$\textsf{Gauss2F1}_q(a-d-n,b-d-n,c-d+1,d+k+2 n+1).$$
	After normalization,
	$$F(n,k) := \frac{q^{(1-a-b+c+d) k+(3-a-b+2 c+2 d) n+2 k n+3 n^2} (q^{1+a};q)_k (q^{1+b};q)_k (q^{1-b+c},q^{1-a+d},q^{1-a+c},q^{1+a+k},q^{1+b+k},q^{1-b+d};q)_n}{(1-q^{1+c+k+2 n}) (1-q^{1+d+k+2 n}) (q^{1+c};q)_k (q^{1+d};q)_k (q^{1-a-b+c+d},q^{1+c+k},q^{1+d+k};q)_{2 n}}.$$
	Now apply Proposition \ref{WZ_prop}, it is easy to see
	$$g(n) = \lim_{k\to \infty} G(n,k) = 0, \qquad \lim_{n\to \infty} \sum_{k\geq 0}F(n,k) = 0,$$
	because of the exponentially decaying factor $q^{2kn+3n^2}$ in the numerator. The equality $\sum_{k\geq 0} F(0,k) = \sum_{n\geq 0} G(n,0)$ yields
	\begin{multline} \label{ex1Q} \sum_{k\geq 0} \frac{q^{k(1-a-b+c+d)} (q^{a+1};q)_k (q^{b+1};q)_k }{(q^{c+1};q)_{k+1} (q^{d+1};q)_{k+1}} \\  = \sum_{n\geq 1} \frac{q^{-1 - a - b - 2 c - 2 d + (-3 - a - b + 2 c + 2 d) n + 3 n^2}    P_q(n) (q^{a+1},q^{b+1},q^{-a+c+1},q^{-a+d+1},q^{-b+c+1},q^{-b+d+1};q)_{n-1}}{ (q^{c+1},q^{d+1},q^{-a-b+c+d+1};q)_{2 n}} \end{multline}
	where {\small $P_q(n) = -q^{a+b+c+d+2 n+1}-q^{2 a+b+c+d+3 n}-q^{a+2 b+c+d+3 n}+q^{2 a+2 b+c+d+4 n}+q^{2 a+2 b+c+d+4 n+1}+q^{a+b+2 c+d+4 n}+q^{a+b+2 c+d+4 n+1}+q^{a+b+c+2 d+4 n}+q^{a+b+c+2 d+4 n+1}-q^{a+b+2 c+2 d+6 n+1}-q^{2 a+b+2 c+2 d+7 n}-q^{a+2 b+2 c+2 d+7 n}+q^{a+b+3 c+2 d+8 n}+q^{a+b+2 c+3 d+8 n}-q^{2 a+2 b+c+2 n+1}-q^{2 a+2 b+d+2 n+1}+q^{2 a+2 b+1}+q^{2 a+2 c+2 d+6 n}-q^{a+3 c+2 d+7 n}-q^{a+2 c+3 d+7 n}+q^{2 b+2 c+2 d+6 n}-q^{b+3 c+2 d+7 n}-q^{b+2 c+3 d+7 n}+q^{3 c+3 d+8 n}$}. This is a $q$-analogue of \eqref{ex1}. Making $a=b=c=d=0$ in \eqref{ex1Q} produces
	$$\sum_{k\geq 0} \frac{q^k}{(1-q^{1+k})^2} = \sum_{n\geq 1} \frac{q^{3 n^2-3 n} (1-2 q^{3 n-1}-3 q^{4 n-1}-3 q^{5 n-1}+3 q^{n}+3 q^{2 n}+q^{3 n}) (q;q)_n^6}{(1-q^n)^3 (q;q)_{2 n}^3}$$
	this is a $q$-analogue of \eqref{ex1S}.
\end{example}

From the preceding examples, we see that to obtain a $q$-analogue of a hypergeometric identity of the form
\begin{equation}\label{to_be_proved_identity_0}
	\sum_{n\ge 0} z^n \frac{(a_1)_n\cdots (a_m)_n}{(b_1)_n \cdots (b_m)_n} R(n)
	= \text{a simple constant},
	\qquad
	z\in\mathbb{Q}, \quad 
	a_i,b_i\in \mathbb{Q}\cap(0,1], \quad
	R(n)\in \mathbb{Q}(n),
\end{equation}
it suffices to find a (classical) WZ-seed that induces \eqref{to_be_proved_identity_0}, and then lift this seed to the $q$-setting.

A natural question arises: given an identity such as \eqref{to_be_proved_identity_0}, how does one find a classical WZ-seed that induces it in the first place? In general, there is no simple answer. At present, the most effective approach is an exhaustive search of a combinatorial nature, as described in \cite[Section~3.4]{au2025wilf}. When this search fails, as it does for certain $1/\pi^k$-formulas, no $q$-analogue can be produced using the current collection of WZ-seeds.

\subsection{Some computational remarks}
The application of the WZ method to infinite series is inherently computationally intensive. The additional conceptual layer introduced by WZ-seeds further increases this complexity. It is therefore worthwhile to explain how the various steps were implemented in order to obtain the examples presented in this paper.

Broadly speaking, the procedure consists of two types of tasks: those that are purely algebraic or combinatorial, and those that require analytic arguments. Recall that our goal is to start from a classical hypergeometric summation formula and construct a multivariable $q$-analogue. The algebraic part consists of the following three steps:
\begin{enumerate}[leftmargin=*]
	\item searching for suitable auxiliary parameters together with a WZ-seed that reproduces the exponential and Pochhammer components of the original formula;
	\item normalizing this WZ-seed into its $q$-version;
	\item computing the corresponding WZ-mate $G(n,k)$.
\end{enumerate}

A systematic procedure for the first step is described in \cite[Section~3.4]{au2025wilf}; the second is straightforward to implement, and the third was described at the end of the previous subsection. These steps are routine and can be largely automated. They have been implemented in \textit{Mathematica}; the corresponding code is available at
\url{https://sites.google.com/view/kc-au/2403-04555}.

The analytic part is less straightforward. Given a WZ-pair $(F(n,k),G(n,k))$, there is currently no general algorithm that automatically produces the desired infinite series. In the examples of this paper, two main approaches are employed:
\begin{itemize}[leftmargin=*]
	\item Proposition \ref{WZ_prop}, which is used in the first twelve examples;
	\item Guillera's method of flawless WZ-pairs \cite{guillera2025wz}, which is used from Example~XIII onward.
\end{itemize}

For either method, analytic complications may arise and often require delicate manipulations or \textit{ad hoc} arguments. This analytic component remains largely case-dependent and, at present, no satisfactory method for automating it is known. Indeed, even in the classical hypergeometric setting, a wide variety of analytic difficulties already occur; see \cite{au2025wilf}.

For readers interested in the computational aspects, we also provide detailed \textit{Mathematica} code for three representative examples (Examples~\ref{ex_1}, \ref{ex_13}, and \ref{ex_17}). These examples illustrate the procedure step by step and provide a more hands-on understanding of the computations involved.

\section{$q$-analogue of Ramanujan-type $1/\pi^k$ formulas: 21 examples}
In the second part of this article, we present $21$ examples, each consisting of a $q$-analogue of a Ramanujan-type $1/\pi^k$ formula. This collection subsumes all previously known results in the literature.

Some of these identities, especially their multi-variable extensions, are rather lengthy. For the convenience of the reader, copy-and-pasteable versions of all formulas, written in \textit{Mathematica} syntax, are provided\footnote{\url{https://sites.google.com/view/kc-au/2403-04555}}. For readers without access to \textit{Mathematica}, we also provide a plain-text repository of all formulas\footnote{\url{https://raw.githubusercontent.com/pisco125/creative_telescoping/refs/heads/main/Repository_of_formulas_tidy.txt}}. While this version sacrifices mathematical typesetting, it offers a convenient and software-independent reference.\par

We remind the reader of our convention that 0<q<1. This applies to all q-series in the sequel.

\begin{romexample}\label{ex_1}
This example concerns a $q$-analogue of \begin{equation}\label{equation_ex_1}\sum_{n\geq 0} \left(-\frac{3^3}{2^9}\right)^n \frac{\left(\frac{1}{2}\right)_n \left(\frac{1}{6}\right)_n \left(\frac{5}{6}\right)_n}{(1)_n^3} (154n+15) = \frac{32\sqrt{2}}{\pi}.\end{equation} Let us first look at the classical hypergeometric case and establish a multi-variable generalization:
	\begin{multline*}
		a\sum_{k\geq 0} \frac{2^{-2 k} \left(-a+b+c+\frac{1}{2}\right)_{2 k}}{(b+1)_k (c+1)_k} = \sum_{n\geq 0} \frac{(-1)^n 2^{-5 n-6} P(n)\left(a+b-c+\frac{1}{2}\right)_n \left(a-b+c+\frac{1}{2}\right)_n \left(-a+b+c+\frac{1}{2}\right)_{3 n}}{(a+1)_n (b+1)_{2 n+1} (c+1)_{2 n+1}} 
	\end{multline*}
	where {\small $P(n) = 8 a^3-8 a^2 b-8 a^2 c-8 a^2 n+4 a^2-8 a b^2+48 a b c+48 a b n+24 a b-8 a c^2+48 a c n+24 a c+88 a n^2+104 a n+30 a+8 b^3+24 b^2 c+88 b^2 n+36 b^2+24 b c^2+240 b c n+72 b c+440 b n^2+312 b n+46 b+8 c^3+88 c^2 n+36 c^2+440 c n^2+312 c n+46 c+616 n^3+676 n^2+214 n+15$}. \par
	How could one (discover and) prove this formula? By searching through our listed WZ-seeds, one is led to consider $$\textsf{Seed1}\left(\frac12- a + b - c - n, 1 + b - c, c + k + 2 n\right),$$ after normalizing, we obtain
	$$F(n,k) := \frac{(-1)^n (a+n) 2^{-2 k-5 n} \left(-a+b+c+\frac{1}{2}\right)_{2 k} \left(a+b-c+\frac{1}{2}\right)_n \left(a-b+c+\frac{1}{2}\right)_n \left(-a+b+c+2 k+\frac{1}{2}\right)_{3 n}}{(a+1)_n (b+1)_k (c+1)_k (b+k+1)_{2 n} (c+k+1)_{2 n}}.$$
	One finds its WZ-mate $G(n,k)$. When applying Proposition \ref{WZ_prop}, we have $g(n) = 0$ and $\lim_{n\to \infty} \sum_{k\geq 0}F(n,k) = 0$; hence
	$$\sum_{k\geq 0} F(0,k) = \sum_{n\geq 0} G(n,0).$$
	This is exactly the above formula. The summand in $k$ decays at rate $O(k^{-1-a})$ at infinity, so we have to assume $a>0$ for convergence of the LHS. Now we want to let $a\to 0$ from the right. To do this, we rewrite the $k$-summand as 
	$$\frac{2^{-2 k} \left(-a+b+c+\frac{1}{2}\right)_{2 k}}{(b+1)_k (c+1)_k} = \frac{2^{-a+b+c-\frac{1}{2}}  \Gamma (b+1) \Gamma (c+1)}{\sqrt{\pi } \Gamma \left(-a+b+c+\frac{1}{2}\right)}k^{-a-1} + O(k^{-2-a}),\quad k\to\infty,$$ 
	therefore the limit of the LHS as $a\to 0^+$ is
	$$\begin{aligned}&\lim_{a\to 0^+} \left( a\sum_{k\geq 0} \frac{2^{-2 k} \left(-a+b+c+\frac{1}{2}\right)_{2 k}}{(b+1)_k (c+1)_k} \right) = \frac{2^{b+c-\frac12} \Gamma (b+1) \Gamma (c+1)}{\sqrt{\pi } \Gamma \left(b+c+\frac{1}{2}\right)} \underbrace{\lim_{a\to 0^+} a \sum_{k\geq 0} k^{-a-1}}_{=1},\end{aligned}$$
	where the last limit is $1$ because $\lim_{a\to 0} a\zeta(a+1) = 1$. Letting $a\to 0$ on the RHS poses no issue, and we obtain the following two-variable identity:
	\begin{equation}\label{aux_1}\sum_{n\geq 0} \frac{(-1)^n 2^{-5 n-6} \tilde{P}(n) \left(b-c+\frac{1}{2}\right)_n \left(-b+c+\frac{1}{2}\right)_n \left(b+c+\frac{1}{2}\right)_{3 n}}{(b+2 n+1) (c+2 n+1) (1)_n (b+1)_{2 n} (c+1)_{2 n}} = \frac{2^{b+c-\frac12} \Gamma (b+1) \Gamma (c+1)}{\sqrt{\pi } \Gamma \left(b+c+\frac{1}{2}\right)}\end{equation}
	with {\small $\tilde{P}(n) =  8 b^3+24 b^2 c+88 b^2 n+36 b^2+24 b c^2+240 b c n+72 b c+440 b n^2+312 b n+46 b+8 c^3+88 c^2 n+36 c^2+440 c n^2+312 c n+46 c+616 n^3+676 n^2+214 n+15$}. 
	Specializing to $b=c=0$ recovers the $1/\pi$-formula \eqref{equation_ex_1}. \par
	
	We are mainly interested in the $q$-version of the above argument, which can be carried out in exact parallel. Start from the $q$-counterpart of the same WZ-seed, namely\footnote{we normalize $q$ by $q^2$ so that the final formula we shall obtain contains only integral power, instead of half-integral, power of $q$.}
	$$\textsf{Seed1}_{q^2}\left(\frac12- a + b - c - n, 1 + b - c, c + k + 2 n\right).$$
	Normalizing gives 
$$F(n,k) := \frac{\splitfrac{(-1)^n(1-q^{2a+4n})q^{2ak+3an+bn+cn+4kn+7n^2}(q^{-a+b+c+1};q^2)_{2k}}{\times(q^{a+b-c+1},q^{a-b+c+1};q^2)_n(q^{-a+b+c+4k+1};q^2)_{3n}}}{(q^{2a+4};q^4)_n(q^{2b+4},q^{2c+4};q^4)_k(q^{2b+4k+4},q^{2c+4k+4};q^4)_{2n}}.$$
One finds its WZ-mate $G(n,k)$. When applying Proposition \ref{WZ_prop}, we have $g(n) = 0$ and $\lim_{n\to \infty} \sum_{k\geq 0}F(n,k) = 0$; both are easy to see because of the $q^{7n^2}$ factor in the numerator. The equality $\sum_{k\geq 0} F(0,k) = \sum_{n\geq 0} G(n,0)$
gives 
\begin{multline}\label{aux_14}(1-q^{2a})\sum_{k\geq0}\frac{q^{2ak}(q^{-a+b+c+1};q^2)_{2k}}{(q^{2b+4};q^4)_k(q^{2c+4};q^4)_k}\\=\sum_{n\geq0}\frac{(-1)^nP(n)q^{3an+bn+cn+7n^2}(q^{a+b-c+1},q^{a-b+c+1};q^2)_n(q^{-a+b+c+1};q^2)_{3n}}{(q^{2b+8n+4}-1)(q^{2c+8n+4}-1)(q^{2a+4},q^{2b+4},q^{2c+4};q^4)_n(q^{2b+4};q^4)_{2n}(q^{2c+4};q^4)_{2n}},\end{multline}
where {\small $P(n) = 1-q^{4+2b+8n}-q^{4+2c+8n}-q^{1+a+b+c+10n}-q^{3+a+b+c+10n}+q^{4+2a+2b+12n}+q^{4+2a+2c+12n}-q^{5+3a+b+c+14n}+q^{4+2b+2c+16n}+q^{8+2b+2c+16n}+q^{6+2a+2b+2c+20n}-q^{9+a+3b+3c+26n}$}. Here the same restriction on $a>0$ is needed to ensure convergence of the LHS. Now we find the limit of LHS as $a\to 0^+$,  since $$\frac{q^{2ak}(q^{-a+b+c+1};q^2)_{2k}}{(q^{2b+4};q^4)_k(q^{2c+4};q^4)_k}=q^{2ak}\frac{(q^{-a+b+c+1};q^2)_\infty}{(q^{2b+4};q^4)_\infty(q^{2c+4};q^4)_\infty}+O(q^{2ak}q^k),\quad k\to \infty,$$
we see 	\begin{multline}\label{aux_15}\lim_{a\to0^+}\left((1-q^{2a})\sum_{k\geq0}\frac{q^{2ak}(q^{-a+b+c+1};q^2)_{2k}}{(q^{2b+4};q^4)_k(q^{2c+4};q^4)_k}\right)=\frac{(q^{b+c+1};q^2)_\infty}{(q^{2b+4};q^4)_\infty(q^{2c+4};q^4)_\infty}\lim_{a\to0^+}(1-q^{2a})\sum_{k\geq0}q^{2ak}\\=\frac{(q^{b+c+1};q^2)_\infty}{(q^{2b+4};q^4)_\infty(q^{2c+4};q^4)_\infty}.\end{multline}
Letting $a=0$ on the RHS poses no issue, so we arrive at
$$\sum_{n\geq0}\frac{(-1)^nq^{bn+cn+7n^2}\tilde{P}(n)(q^{1+b-c},q^{1-b+c};q^2)_n(q^{1+b+c};q^2)_{3n}}{(1-q^{4+2b+8n})(1-q^{4+2c+8n})(q^4,q^{4+2b},q^{4+2c};q^4)_n}=\frac{(q^{b+c+1};q^2)_\infty}{(q^{2b+4};q^4)_\infty(q^{2c+4};q^4)_\infty},$$
where {\small $\tilde{P}(n)=1-q^{4+2b+8n}-q^{4+2c+8n}-q^{1+b+c+10n}-q^{3+b+c+10n}+q^{4+2b+12n}+q^{4+2c+12n}-q^{5+b+c+14n}+q^{4+2b+2c+16n}+q^{8+2b+2c+16n}+q^{6+2b+2c+20n}-q^{9+3b+3c+26n}$}.  
	This is a $q$-analogue of identity \eqref{aux_1}. Setting $b=c=0$ produces
$$\sum_{n\geq 0} \frac{\splitfrac{(-1)^{n} q^{7 n^2}  \big(1+2 q^{1+2 n}+3 q^{2+4 n}+4 q^{3+6 n}+3 q^{4+8 n}-q^{1+10 n}-q^{3+10 n}+2 q^{5+10 n}}{-2 q^{2+12 n}+q^{6+12 n}-3 q^{3+14 n}-3 q^{4+16 n}-3 q^{5+18 n}-2 q^{6+20 n}-q^{7+22 n}\big)(q;q^2)_{3 n} (q;q^2)_n^2 }}{(q^{2 n+1}+1)^2 (q^{4 n+2}+1)^2 (q^4;q^4)_n (q^4;q^4)_{2 n}^2} = \frac{(q;q^2)_\infty}{(q^4;q^4)_\infty^2}.
$$
This is a $q$-analogue of the $1/\pi$-formula \eqref{equation_ex_1}. A detailed workflow of this example is available in the accompanying \textit{Mathematica} file.
\end{romexample}

A one-parameter special case of the above two-parameter equation was derived in \cite{guillera2018wz}.

In all previous examples, we juxtaposed the classical and $q$-case for pedagogical purposes. From now on, we shall focus exclusively on the $q$-case, which is a more general setting since the classical formulas are recovered under $q\to 1$. \par

The next eight examples all follow the same template. Starting from a given WZ-seed involving accessory variables $a,b,c,d,\cdots$, one first derives an identity of the form \eqref{aux_14}. One then lets the parameter $a\to 0^+$, which yields a summation formula depending on the remaining free variables $b,c,d,\cdots$. To avoid unnecessary repetition, in Examples~\ref{ex_2}–\ref{ex_9} we present only the following data:
\begin{itemize}[leftmargin=*]
	\item the underlying WZ-seed;
	\item the resulting summation formula with free parameters $b,c,d$;
	\item the specialization of this formula obtained by setting all parameters equal to $0$;
	\item the corresponding classical $1/\pi^k$-formula, recovered in the limit $q\to 1$.
\end{itemize}

\begin{romexample}\label{ex_2}
Starting from $\textsf{Gauss2F1}_{q^2}(-a-b+c-n+\frac{1}{2},-a+c-n+\frac{1}{2},-a-b+2 c-n+1,2 a+2 b-2 c+d+k+2 n)$, we get the formula
\begin{multline*}\sum_{n\geq0}\frac{\splitfrac{q^{bn+dn+2n^2}(q^{2c}-q^{2+2b+d+4n}-q^{1+b+c+d+4n}-q^{1+2b+c+d+4n}+q^{2+2b+d+6n}+q^{2+3b+2d+6n})}{\times(q^{1+2b-c+d},q^{1-c},q^{1+b-c},q^{1+b-c+d};q^2)_n}}{\left(q^{2c}-q^{2+2b+d+4n}\right)(q^2,q^{2+b+d};q^2)_n(q^{2+2b-2c+d};q^2)_{2n}}\\=\frac{(q^{b-c+d+1},q^{2b-c+d+1};q^2)_{\infty}}{(q^{b+d+2},q^{2b-2c+d+2};q^2)_{\infty}}.\end{multline*} \\
Specializing to $b=c=d=0$: $$ \sum_{n\geq 0} \frac{q^{2 n^2} (1+q^{2 n+1}-2 q^{4 n+1}) (q;q^2)_n^4}{(q^{2 n+1}+1) (q^2;q^2)_{2 n} (q^2;q^2)_n^2} = \frac{(q;q^2)_{\infty }^2}{(q^2;q^2)_{\infty }^2}.$$
This is a $q$-analogue of $$\sum_{n\geq 0} \left(\frac{1}{2^2}\right)^n \frac{\left(\frac{1}{2}\right)_n^3}{(1)_n^3} (6n+1) = \frac{4}{\pi}.$$
\end{romexample}

\begin{romexample}\label{ex_3}
Starting from $\textsf{Gauss2F1}_{q^2}(-a-b+c-n+\frac{1}{2},b-n+\frac{1}{2},c+1,d+k+2 n)$, we get the formula
\begin{multline*}\sum_{n\geq0}\frac{q^{cn+2dn+6n^2}\tilde{P}(n)(q^{1-b},q^{1+b},q^{1+b-c},q^{1-b+c},q^{1-b+c+d},q^{1+b+d};q^2)_n}{(1-q^{2+4n})(1-q^{2+d+4n})(1-q^{2+c+d+4n})(q^2,q^{2+d},q^{2+c+d};q^2)_{2n}}=\frac{(q^{-b+c+d+1},q^{b+d+1};q^2)_{\infty}}{(q^{c+d+2},q^{d+2};q^2)_\infty},\end{multline*}
where {\small $\tilde{P}(n)=-(q^{2b}+q^{c})q^{-b+d+6n+1}+(-q^{2b+c+2}+q^{4b}+q^{2c})q^{-2b+2d+12n+4}-(q^{2b}+q^{c})(q^{c+d}+q^{c}+1)q^{-b+2d+14n+5}-q^{4n+2}(q^{c+d}+q^{d}+1)+(q^{c+d}+q^{d}+1)q^{c+2d+16n+6}+(q^2+1)(q^{c+d}+q^{c}+1)q^{d+8n+2}+1$}. \\
Specializing to $b=c=d=0$: $$\sum_{n\geq 0}\frac{q^{6 n^2} (3 q^{2 n+1}+3 q^{4 n+2}-2 q^{6 n+1}+q^{6 n+3}-3 q^{8 n+2}-3 q^{10 n+3}+1) (q;q^2)_n^6}{(q^{2 n+1}+1)^3 (q^2;q^2)_{2 n}^3} = \frac{(q;q^2)_{\infty }^2}{(q^2;q^2)_{\infty }^2}.$$
This is a $q$-analogue of $$\sum_{n\geq 0} \left(\frac{1}{2^6}\right)^n \frac{\left(\frac{1}{2}\right)_n^3}{(1)_n^3} (42n+5) = \frac{16}{\pi}.$$
This formula can also be found in \cite{chen2021hidden}.
\end{romexample}

\begin{romexample}\label{ex_4}
	Using $\textsf{Seed1}_{q^2}(-a+b-n+\frac{1}{2},b+1,c+k+n)$ gives
	$$\sum_{n\geq0}\frac{(-1)^nq^{bn+2cn+3n^2}(1-q^{b+2c+6n+1})(q^{1-b},q^{b+1},q^{b+2c+1};q^2)_n}{(q^4,q^{2c+4},q^{2b+2c+4};q^4)_n}=\frac{(q^{b+2c+1};q^2)_{\infty}}{(q^{2c+4};q^4)_{\infty}(q^{2b+2c+4};q^4)_{\infty}}.$$ \\
Specializing to $b=c=0$: $$\sum_{n\geq 0}\frac{(-1)^n q^{3 n^2} (1-q^{6 n+1}) (q;q^2)_n^3}{(q^4;q^4)_n^3} = \frac{(q;q^2)_{\infty }}{(q^4;q^4)_\infty^2} .$$
	This is a $q$-analogue of $$\sum_{n\geq 0} \left(\frac{-1}{2^3}\right)^n \frac{\left(\frac{1}{2}\right)_n^3}{(1)_n^3} (6n+1) = \frac{2\sqrt{2}}{\pi}.$$
A one-parameter special case of the above example can be found in \cite{guillera2018wz, guo2018q}.
\end{romexample}

\begin{romexample}\label{ex_5}
	Using $\textsf{Gauss2F1}_{q^2}(-a-b+c+\frac{1}{2},b-n+\frac{1}{2},c+n+1,d+k+n)$ gives
	$$\sum_{n\geq0}\frac{(-1)^nq^{bn+2dn+3n^2}\tilde{P}(n)(q^{1-b},q^{1-b+c+d},q^{1+b};q^2)_n(q^{1-b+c};q^2)_{2n}}{(1-q^{2+4n})(1-q^{2+c+d+4n})(q^2,q^{2+c+d};q^2)_{2n}(q^{2+d};q^2)_n}=\frac{(q^{b+d+1},q^{-b+c+d+1};q^2)_\infty}{(q^{d+2},q^{c+d+2};q^2)_\infty}.$$ \\
	where {\small $\tilde{P}(n) = 1+q^{4-2b+2c+2d+12n}-q^{5-b+2c+2d+14n}+q^{2+c+d+8n}(1+q^2-q^{1-b+d})+q^{1+d+6n}(q-q^{-b+c}+q^{1+c+d})-q^{1+4n}(q+q^{b+d}+q^{1+c+d})$}.\\
	Specializing to $b=c=d=0$: $$\sum_{n\geq 0}\frac{(-1)^n q^{3 n^2} (2 q^{2 n+1}-q^{4 n+1}+q^{4 n+2}-q^{6 n+1}-q^{8 n+2}-q^{10 n+3}+1) (q;q^2)_{2 n} (q;q^2)_n^3}{(q^{2 n+1}+1)^2 (q^2;q^2){}_n (q^2;q^2)_{2 n}^2} = \frac{(q;q^2)_\infty^2}{(q^2;q^2)_\infty^2} .$$
	This is a $q$-analogue of $$\sum_{n\geq 0} \left(\frac{-1}{2^2}\right)^n \frac{\left(\frac{1}{2}\right)_n \left(\frac14\right)_n \left(\frac34\right)_n}{(1)_n^3} (20n+3) = \frac{8}{\pi}.$$
\end{romexample}

\begin{romexample}\label{ex_6}
	Using $\textsf{Dougall5F4}_{q^4}(a+b+c+d,b+\frac{1}{2},c+\frac{1}{4},d-n+\frac{1}{4},e+k+n+\frac{1}{4})$ gives
	\begin{multline*}\sum_{n\geq0}\frac{(-1)^nq^{2dn+2en+2n^2}\tilde{P}(n)(q^{2+2b},q^{1+2c},q^{3-2d},q^{3+2b+2e},q^{2+2c+2e},q^{1+2b+2c+2d+2e};q^4)_n}{(1-q^{1+2e+4n})(1-q^{4+2b+2c+2e+8n})(q^4,q^{1+2e},q^{4+2b+2d+2e},q^{3+2c+2d+2e};q^4)_n(q^{4+2b+2c+2e};q^4)_{2n}}\\=\frac{(q^{3+2b+2e},q^{2+2c+2e},q^{2+2d+2e},q^{1+2b+2c+2d+2e};q^4)_{\infty}}{(q^{1+2e},q^{4+2b+2c+2e},q^{4+2b+2d+2e},q^{3+2c+2d+2e};q^4)_{\infty}},\end{multline*} \\
	where {\small $\tilde{P}(n)=1-q^{1+2e+8n}(q^{2+2b}+q^{1+2c}+q^{3+2b+2c}+q^{2b+2c+2d}+q^{1+2b+2c+2d+2e})+q^{2+2b+2c+2e+12n}(q^2+q^{3+2e}+q^{2d+2e}+q^{2+2b+2d+2e}+q^{1+2c+2d+2e})-q^{6+4b+4c+2d+6e+20n}$}.\\
	Specializing to $b=c=d=e=0$: $$\sum_{n\geq 0}\frac{(-1)^n q^{2 n^2} (q^{4 n+1}+q^{4 n+2}-q^{8 n+1}-q^{8 n+2}-q^{12 n+3}+1) (q;q^2)_{2 n} (q^2;q^4)_n^2}{(q^{4 n+2}+1) (q^4;q^4)_{2 n} (q^4;q^4)_n^2}  = \frac{(q^2;q^4)_\infty^2}{(q^4;q^4)_\infty^2} .$$
	This is another $q$-analogue, different from that in the previous example, of $$\sum_{n\geq 0} \left(\frac{-1}{2^2}\right)^n \frac{\left(\frac{1}{2}\right)_n \left(\frac14\right)_n \left(\frac34\right)_n}{(1)_n^3} (20n+3) = \frac{8}{\pi}.$$

A one-parameter special case of the above example can be found in \cite{guo2020q}. 
\end{romexample}

\begin{romexample}\label{ex_7}
	Using $\textsf{Dougall5F4}_{q^2}(a+b+c+d-n+\frac{1}{2},b-n+\frac{1}{2},c-n+\frac{1}{2},d-n+\frac{1}{2},e+k+n)$ gives
	\begin{multline*}
	\sum_{n\geq0}\frac{(-1)^nq^{bn+cn+dn+2en+n^2}\tilde{P}(n)(q^{1-b},q^{1+b},q^{1-c},q^{1+c},q^{1-d},q^{1+d};q^2)_n}{(1-q^{2+4n})(q^2;q^2)_{2n}(q^{2+e},q^{2+b+c+e},q^{2+c+d+e},q^{2+b+d+e};q^2)_n}\\=\frac{(q^{1+b+e},q^{1+c+e},q^{1+d+e},q^{1+b+c+d+e};q^2)_{\infty}}{(q^{2+e},q^{2+b+c+e},q^{2+c+d+e},q^{2+b+d+e};q^2)_\infty},
	\end{multline*}
	where {\small $\tilde{P}(n)=1-q^{1+4n}(q+q^{b+e}+q^{c+e}+q^{d+e}+q^{b+c+d+e})+q^{1+e+6n}(q+q^{1+b+c}+q^{1+b+d}+q^{1+c+d}+q^{b+c+d+e})-q^{3+b+c+d+2e+10n}$}. \\
	Specializing to $b=c=d=e=0$: 
	$$\sum_{n\geq 0} \frac{(-1)^n q^{n^2} (q^{2 n+1}-4 q^{4 n+1}+q^{6 n+1}+q^{8 n+2}+1) (q;q^2)_n^6}{(q^{2 n+1}+1) (q^2;q^2)_{2 n} (q^2;q^2)_n^4} = \frac{(q;q^2)_\infty^4}{(q^2;q^2)_\infty^4}.$$
		This is a $q$-analogue of $$\sum_{n\geq 0} \left(\frac{-1}{2^2}\right)^n \frac{\left(\frac{1}{2}\right)_n^5}{(1)_n^5}  (20n^2+8n+1) = \frac{8}{\pi^2} .$$
\end{romexample}

\begin{romexample}\label{ex_8}
	Using $\textsf{Dougall5F4}_{q^2}(a+b+c+d-n+\frac{1}{2},b-n+\frac{1}{2},c-n+\frac{1}{2},d-n+\frac{1}{2},e+k+2 n)$ gives
	\begin{multline*}
		\sum_{n\geq0}\frac{(-1)^nq^{2(b+c+d+e)n+5n^2} \tilde{P}(n) (q^{1-b},q^{1+b},q^{1-c},q^{1+c},q^{1-d},q^{1+b+e},q^{1+c+e},q^{1+d+e},q^{1+d},q^{1+b+c+d+e};q^2)_n}{(q^2,q^{2+e},q^{2+b+c+e},q^{2+c+d+e},q^{2+b+d+e};q^2)_{2n+2}}\\=\frac{(q^{1+b+e},q^{1+c+e},q^{1+d+e},q^{1+b+c+d+e};q^2)_{\infty}}{(q^{2+e},q^{2+b+c+e},q^{2+c+d+e},q^{2+b+d+e};q^2)_\infty},
	\end{multline*}
	where $\tilde{P}(n)\in \mathbb{Z}[q,q^b,q^c,q^d,q^e]$ is a long $q$-polynomial\footnote{whose full expression can be found in the accompanying \textit{Mathematica} file}. Specializing to $b=c=d=e=0$: 
	$$\sum_{n\geq 0} \frac{\splitfrac{(-1)^n q^{5 n^2} \big[ 1+5 q^{1+2 n}+10 q^{2+4 n}+5 q^{1+6 n} (-1+2 q^2)+5 q^{2+8 n} (-3+q^2)+q^{10 n} (q-24 q^3+q^5)}{+q^{12 n} (5 q^2-15 q^4)-5 q^{3+14 n} (-2+q^2)+10 q^{4+16 n}+5 q^{5+18 n}+q^{6+20 n} \big]  (q;q^2)_n^{10}}}{(q^{2 n+1}+1)^5 (q^2;q^2)_{2 n}^5}=  \frac{(q;q^2)_\infty^4}{(q^2;q^2)_\infty^4}.$$
	This is a $q$-analogue of $$\sum_{n\geq 0} \left(\frac{-1}{2^{10}}\right)^n \frac{\left(\frac{1}{2}\right)_n^5}{(1)_n^5} (820n^2+180n+13) = \frac{128}{\pi^2}  .$$
\end{romexample}

\begin{romexample}\label{ex_9}
	Using $\textsf{Dougall5F4}_{q^2}(a+b+c+d-n+\frac{1}{2},b-n+\frac{1}{2},c-n+\frac{1}{2},d+\frac{1}{2},e+k+2 n)$ gives
	\begin{multline*}
		\sum_{n\geq0}\frac{q^{bn+cn+en+2n^2}\tilde{P}(n)(q^{1-b},q^{1-c},q^{1+d},q^{1+b+c+d+e},q^{1+d+e},q^{1+b+e},q^{1+c+e};q^2)_n(q^{1+d+e};q^2)_{2n}}{(q^2,q^{2+b+c+e};q^2)_n(q^{2+e},q^{2+c+d+e},q^{2+b+d+e};q^2)_{2n+2}}\\=\frac{(q^{1+b+e},q^{1+c+e},q^{1+d+e},q^{1+b+c+d+e};q^2)_{\infty}}{(q^{2+e},q^{2+b+c+e},q^{2+c+d+e},q^{2+b+d+e};q^2)_\infty},
	\end{multline*}
	where {\small $\tilde{P}(n)=1-q^{1+e+4n}(q+q^{b}+q^{c}+q^{1+b+d}+q^{1+c+d}+q^{b+c+d})+q^{1+e+6n}(q-q^{d}+q^{1+b+d}+q^{1+c+d})(1+q^{b+c+e})+q^{1+d+e+8n}(-q^2+q^{1+b+e}+q^{3+b+e}+q^{1+c+e}+q^{3+c+e}+q^{b+c+e}+q^{1+b+c+d+e}+q^{3+b+c+d+e}-q^{2+2b+2c+2e})+q^{3+d+2e+12n}(q^{1+d}-q^{b+c+e}-q^{2+b+c+e}-q^{3+b+c+d+e}-q^{2b+c+d+e}-q^{2+2b+c+d+e}-q^{b+2c+d+e}-q^{2+b+2c+d+e}+q^{1+2b+2c+d+2e})-q^{5+2d+3e+14n}(q^{b}+q^{c}-q^{1+b+c}+q^{b+c+d})(1+q^{b+c+e})+q^{5+b+c+2d+4e+16n}(q+q^{b}+q^{c}+q^{1+b+d}+q^{1+c+d}+q^{b+c+d})-q^{7+2b+2c+3d+5e+20n}$}.\\ Specializing to $b=c=d=e=0$: 
	$$\sum_{n\geq 0} \frac{\splitfrac{q^{2 n^2} \big[1+3 q^{1+2 n}+3 (-1+q) q^{1+4 n}+q^{1+6 n} (-2-3 q+q^2)+q^{8 n} (q-3 q^2-2 q^3)}{-3 (-1+q) q^{2+10 n}+3 q^{3+12 n}+q^{4+14 n} \big] (q;q^2)_{2 n} (q;q^2)_n^6}}{(q^{2 n+1}+1)^3 (q^2;q^2)_n^2 (q^2;q^2)_{2 n}^3} =  \frac{(q;q^2)_\infty^4}{(q^2;q^2)_\infty^4}.$$
	This is a $q$-analogue of $$\sum_{n\geq 0} \left(\frac{1}{2^4}\right)^n \frac{\left(\frac{1}{2}\right)_n^3 \left(\frac{1}{4}\right)_n \left(\frac{3}{4}\right)_n}{(1)_n^5} (120n^2 + 34n+3) = \frac{32}{\pi^2}.$$
\end{romexample}
The $q$-analogue of the above two $1/\pi^2$-formulas can also be found in \cite{chu2018q, campbellhal}.

\begin{romexample}\label{ex_10}
We claim the $1/\pi$-formula $$\sum_{n\geq 0} \left(\frac{-1}{48}\right)^n \frac{\left(\frac{1}{2}\right)_n \left(\frac{1}{4}\right)_n \left(\frac{3}{4}\right)_n}{(1)_n^3} (28n+3) = \frac{16}{\sqrt{3}\pi}$$ has a $q$-analogue:
\begin{multline*}\sum_{n\geq 0} \frac{\splitfrac{(-1)^n q^{5 n^2} \big(1+q^{1+2 n}+q^{1+4 n} (1+q)+q^{2+6 n} (1+q)+q^{2+8 n} (1+q)+q^{1+10 n} (q^3-1)}{-q^{2+12 n} (1+q)-q^{2+14 n} (1+q)-q^{3+16 n} (1+q)-q^{4+18 n}-q^{5+20 n}\big) (q;q^2)_n (q^3;q^6)_n (q;q^2)_{2 n}^2}}{(q^{2 n+1}+1) (q^{4 n+1}+1) (q^{4 n+2}+1) (q^2;q^2)_{4 n} (q^6;q^6)_n^2} \\ = \frac{(q;q^2)_{\infty } (q^3;q^6)_{\infty }}{(q^2;q^2)_{\infty } (q^6;q^6)_{\infty}}.\end{multline*}
This in turn has a two-variable extension:
\begin{multline}\label{aux_3}\sum_{n\geq0}\frac{(-1)^nq^{bn+cn+5n^2}P(n)(q^{-b-c+3};q^6)_n(q^{b+c+1};q^2)_n(q^{b+1},q^{1-b};q^2)_{2n}}{(-q^{8n+2}-q^{8n+4}+q^{16n+6}+1)(q^2;q^2)_{4n}(q^{c+6},q^{2b+c+6};q^6)_n}\\=\frac{(q^{3-b},q^{3+b},q^6;q^6)_{\infty}(q^{1+b+c};q^2)_{\infty}}{(q^2;q^2)_{\infty}(q^{6+c},q^{3+b+c},q^{6+2b+c};q^6)_{\infty}},\end{multline}
where {\small $P(n)=1-q^{4+8n}-q^{1+b+c+10n}(1+q^2)-q^{3-b+12n}(1+q^{2b})+q^{4+c+14n}(1+q^{2b})+q^{16n}(q^4+q^6)+q^{3+b+c+18n}-q^{7+b+c+26n}$}. To establish this formula, start from the WZ-seed $\textsf{Seed10}_{q^2}(3 b-2 n+\frac{1}{2},2 b+1,c+k+n)$, after some rewriting, denote
$$F(n,k) := \frac{\splitfrac{(-1)^n(1-q^{8n-2})(1-q^{8n})q^{bn+cn+6kn+5n^2}(q^{b+1},q^{-b+1};q^2)_{2n}}{(q^{b+c+1};q^2)_{3k}(q^{-b-c-6k+3};q^6)_n(q^{b+c+6k+1};q^2)_n}}{(q^{2};q^2)_{4n}(q^{2b+c+6},q^{c+6},q^{c+3};q^6)_k(1-q^{-b-c-6k+6n-3})(q^{c+6k+6},q^{2b+c+6k+6};q^6)_n}.$$
Note that $F(n,0)=0$ and $\lim_{n\to \infty} \sum_{k\geq 0}F(n,k) = 0$ because of $q^{6kn+5n^2}$ in the numerator. The limit $g(n) := \lim_{k\to \infty} G(n,k) $ is however non-zero, as indicated by the term $(q^{-b-c-6 k+3};q^6)_n$. Proposition \ref{WZ_prop} implies $\sum_{n\geq 0} G(n,0) = \sum_{n\geq 0} g(n),$ which is
\begin{multline*}\sum_{n\geq 0} \frac{(-1)^nq^{bn+cn+5n^2}P(n)(q^{-b-c+3};q^6)_n(q^{b+c+1};q^2)_n(q^{b+1},q^{1-b};q^2)_{2n}}{(-q^{8n+2}-q^{8n+4}+q^{16n+6}+1)(q^2;q^2)_{4n}(q^{c+6},q^{2b+c+6};q^6)_n}\\=\frac{(q^{b+c+1};q^2)_{\infty}}{(q^{c+6},q^{b+c+3},q^{2b+c+6};q^6)_{\infty}}\sum_{n\geq0}\frac{q^{8n^2}(-q^{b+12n+3}-q^{-b+12n+3}-q^{8n+4}+q^{16n+4}+q^{16n+6}+1)(q^{1-b},q^{b+1};q^2)_{2n}}{(-q^{4n+2}-q^{8n+2}+q^{12n+4}+1)(q^2;q^2)_{4n}}. \end{multline*}
It remains to evaluate the summation on the RHS: it is independent of $c$. So if we can evaluate the LHS for a special value of $c$, then we will obtain its value. We pick $c=3-b$, so the summation on the LHS terminates, with $n=0$ as the only contributing term because of the term $(q^{-b-c+3};q^6)_n$, thus proving formula \eqref{aux_3}.
\end{romexample}

\begin{romexample}\label{ex_11}
The $1/\pi$-formula $\sum_{n\geq 0} \left(\frac{1}{2^2}\right)^n \frac{\left(\frac{1}{2}\right)_n^3}{(1)_n^3} (6n+1) = \frac{4}{\pi}$
has another $q$-analogue (different than that of Example \ref{ex_2}):
$$\sum_{n\geq 0}\frac{q^{n^2} (1-q^{1+6n}) (q;q^2)_n^2 (q^2;q^4)_n}{(q^4;q^4)_n^3} = \frac{(q^2;q^4)_{\infty }^2}{(q^4;q^4)_{\infty }^2},$$
which follows from the three-variable extension below:
\begin{multline}\label{aux_33} \sum_{n\geq0}\frac{q^{-2bn+cn+n^2}(1-q^{c+2d+6n+1})(q^{2b+2d+2};q^4)_n(q^{1-c},q^{c+1},q^{c+2d+1};q^2)_n}{(q^{-2b+c+1};q^2)_n(q^4,q^{2d+4},q^{2c+2d+4};q^4)_n}\\=\frac{(q^{c+2d+1};q^2)_{\infty}(q^{-2b+2c+2},q^{2-2b};q^4)_{\infty}}{(q^{-2b+c+1};q^2)_{\infty}(q^{2c+2d+4},q^{2d+4};q^4)_{\infty}}.\end{multline}
This can be proved exactly as in the previous example by taking $F(n,k)$ to be  $\textsf{Seed2}_q(c-n+\frac{1}{2},b+\frac{1}{2},c+1,d+k+n)$. The above formula is $\sum_{n\geq 0} G(n,0) = \sum_{n\geq 0} g(n),$ and the $\sum_{n\geq 0} g(n)$ can be found using the same trick as above.
\end{romexample}

A single-variable special case of the above example was proved in \cite{guo2018ramanujan}.

\begin{remark}
Setting $(b,c,d)=(-\frac16,0,0)$ in the formula \eqref{aux_33} gives
$$\sum_{n\geq 0} \frac{q^{n^2+\frac{4 n}{3}} (1-q^{6 n+1}) (q^{2/3};q^4)_n (q;q^2)_n^3}{(1-q) (q^{7/3};q^2)_n (q^4;q^4)_n^3} = \frac{(q+1) \sqrt[3]{q^2+1} \Gamma _{q^2}(\frac{7}{6})}{\Gamma _{q^2}(\frac{1}{2}) \Gamma _{q^4}(\frac{5}{6})^2}.$$
This is a $q$-analogue of $$\sum_{n\geq 0} \left(\frac{1}{2^2}\right)^n \frac{\left(\frac12\right)_n^3}{(1)_n^3} = \frac{\sqrt{3} \Gamma(\frac{1}{3})^6}{2^{8/3} \pi ^4}.$$
\end{remark}

\begin{romexample}\label{ex_12}
	The $1/\pi$-formula $\sum_{n\geq 0} \left(\frac{1}{2^6}\right)^n \frac{\left(\frac{1}{2}\right)_n^3}{(1)_n^3} (42n+5)= \frac{16}{\pi} $
	has another $q$-analogue (different than that of Example \ref{ex_3}):
$$\begin{aligned}
	& \sum_{n\geq 0} \frac{q^{5 n^2} \left(1+2 q^{4 n+2}+q^{6 n+1}-2 q^{8 n+2}+q^{8 n+4}-q^{10 n+1}+2 q^{10 n+3}-q^{12 n+4}-2 q^{14 n+3}-q^{18 n+5}\right) (q;q^2)_n^2 (q^2;q^4)_n^3}{(q^{4 n+2}+1)^2 (q^4;q^4)_n (q^4;q^4)_{2 n}^2} \\ 
	&= \frac{(q^2;q^4)_{\infty }^2}{(q^4;q^4)_{\infty }^2},
\end{aligned}$$
	which follows from the three-variable extension below:
	\begin{multline*}\sum_{n\geq0}\frac{q^{-2bn+cn+5n^2}P(n)(q^{2-2b},q^{2-2b+2c},q^{2+2b+2d};q^4)_n(q^{1-c},q^{1+c};q^2)_n(q^{1+c+2d};q^2)_{3n}}{(q^4;q^4)_n(q^{1-2b+c};q^2)_{3n+2}(q^{4+2d},q^{4+2c+2d};q^4)_{2n+1}}\\=\frac{(q^{c+2d+1};q^2)_{\infty}(q^{-2b+2c+2},q^{2-2b};q^4)_{\infty}}{(q^{-2b+c+1};q^2)_{\infty}(q^{2c+2d+4},q^{2d+4};q^4)_{\infty}},\end{multline*}
	where {\small $P(n)=1-q^{3-2b+c+6n}-q^{7-4b+3c+2d+18n}+q^{6+2c+4d+20n}+q^{10-2b+4c+6d+32n}-q^{13-4b+5c+6d+38n}+q^{3-2b+c+2d+14n}(1+q^2+q^4-q^{2+2b}+q^{2c}+q^{2+2c}+q^{4+2c})-q^{6-4b+2c+4d+24n}(q^{2b}+q^{2+2b}+q^{4+2b}-q^{2+2c}+q^{2b+2c}+q^{2+2b+2c}+q^{4+2b+2c})+q^{1+c+10n}(q^{2-2b}-q^{2d}-q^{2+2d})+q^{2+2c+2d+16n}(q^{-2b}+q^{2+2d}+q^{6+2d})+q^{10-4b+4c+4d+28n}(1+q^2-q^{2b+2d})-q^{2-2b+8n}(1+q^{2c})(1+q^{2+2b+2d})+q^{9-4b+3c+4d+30n}(1+q^{2c})(1+q^{2+2b+2d})-q^{5-4b+3c+2d+22n}(1+q^4+q^{6+2b+2d})-q^{9-4b+3c+4d+26n}(1-q^{2+2b}+q^{2c}+q^{4b+2d})+q^{12n}(q^{4-4b+2c}+q^{4+2d}+q^{4+2c+2d}-q^{2-2b+2c+2d})$}. \par
	This can be proved exactly as in the previous two examples by taking $F(n,k)$ to be  $\textsf{Seed2}_q(c-n+\frac{1}{2},b-n+\frac{1}{2},c+1,d+k+2 n)$. The above formula is $\sum_{n\geq 0} G(n,0) = \sum_{n\geq 0} g(n),$ and the $\sum_{n\geq 0} g(n)$ can be found using the same trick as above.
\end{romexample}

\begin{romexample}\label{ex_13}
The remarkable $1/\pi^4$-formula, conjectured by Cullen and proved by the author in \cite{au2025wilf}
$$\sum _{n\geq 0} \left(\frac{1}{2^{12}}\right)^n \frac{\left(\frac{1}{2}\right)_n^7 \left(\frac{1}{4}\right)_n \left(\frac{3}{4}\right)_n}{(1)_n^9}(43680 n^4+20632 n^3+4340 n^2+466 n+21)=\frac{2048}{\pi ^4}$$
has a $q$-analogue of the form
\begin{equation}\label{aux_4}\sum_{n\geq 0} \frac{q^{4 n^2}P(n) (q^2;q^2)_{4 n} (q;q^2)_n^{16}  }{(q^{4 n}+1) (q^{2 n+1}+1)^9 (q^2;q^2)_{2 n}^{10}} = \frac{(q;q^2)_{\infty }^8}{(q^2;q^2)_{\infty }^8} \end{equation}
where {\small $P(n) =1+ 9 q^{2 n+1} + (36 q^2+1) q^{4 n} + 7(12 q^2-1) q^{6 n+1} + 63 (2 q^2-1) q^{8 n+2} + 21 (6 q^2-11) q^{10 n+3} + 7 (12 q^4-64 q^2+5) q^{12 n+2} + (36 q^4-560 q^2+187) q^{14 n+3} + (9 q^6-441 q^4+509 q^2-9) q^{16 n+2} + (q^6-225 q^4+805 q^2-65) q^{18 n+3} - (65 q^6-805 q^4+225 q^2-1) q^{20 n+2} - (9 q^6-509 q^4+441 q^2-9) q^{22 n+3} + (187 q^4-560 q^2+36) q^{24 n+4} + 7 (5 q^4-64 q^2+12) q^{26 n+5} -21 (11 q^2-6) q^{28 n+6}-63 (q^2-2) q^{30 n+7}-7 (q^2-12) q^{32 n+8} + (q^2+36) q^{34 n+9} + 9 q^{36 n+10} + q^{38 n+11} $}. 

To begin with, one must first identify which WZ-seed proves the classical identity, and then consider its $q$-version, which is
$$\textsf{Dougall7F6}_{q^2}\left(\frac{1}{2}-n,\frac{1}{2}-n,\frac{1}{2}-n,\frac{1}{2}-n,\frac{1}{2}-n,k+2 n\right),$$
so we are inspired to take
$$F(n,k) := \frac{(1-q^{4 n})^5 q^{2 k n+4 n^2} (1-q^{4 k+6 n+1}) (q;q^2)_k^4 (q;q^2)_n^{10} (q^{1-2 k};q^2)_n (q^{2 k+2};q^2)_{4 n} (q^{2 k+1};q^2)_n^5}{(1-q^{-2 k+2 n-1}) (1-q^{2 k+8 n}) (q^2;q^2)_k^4 (q^2;q^2)_{2 n}^5 (q^{2 k+2};q^2)_{2 n}^5}.$$
Then $$G(n,k) = \frac{q^{2 k n+4 n^2} P(n,k) (q;q^2)_k^4 (q;q^2)_n^{10} (q^{1-2 k};q^2)_n (q^{2 k+2};q^2)_{4 n} (q^{2 k+1};q^2)_n^5}{(q^{2 n+1}+1)^5 (q^{2 k+4 n+2}-1)^5 (q^{2 k+8 n}-1)  (q^2;q^2)_k^4 (q^2;q^2)_{2 n}^5 (q^{2 k+2};q^2)_{2 n}^5},$$
where $P(n,k) \in \mathbb{Z}[q,q^n,q^k]$ is a long polynomial which we do not display here. Next we use an analytic procedure due to Guillera (who calls it \textit{flawless WZ-pairs}, \cite{guillera2025wz}) that we will repeat several times for later examples. For fixed $k\in \mathbb{C}$, we have $F(0,k) =0$ and $\lim_{n\to\infty} F(n,k) = 0$. Therefore $$\sum_{n\geq 0} (G(n,k+1) - G(n,k)) = \sum_{n\geq 0} (F(n+1,k) - F(n,k)) = \lim_{n\to\infty} F(n,k) - F(0,k) = 0.$$
That is, $S(k) := \sum_{n\geq 0} G(n,k)$ is 1-periodic in $k$. We claim that $S(k)$ is holomorphic in $k$: this follows from the appearance of $G(n,k)$ if $\operatorname{Re}(k)$ is sufficiently large; periodicity extends this to all $k$. Also, as a function of $q^k$, $S(k)$ has another period: $S(k+ \frac{2\pi i}{\log q}) = S(k)$. Therefore $S$ is a holomorphic double-periodic function, so is a constant. Its value can be obtained by taking $k=1/2$, in which the sum terminates giving only a single contributing $n$: $$S(\frac12) = \sum_{n\geq 0} G(n,\frac12) = G(0,\frac12).$$
Our $q$-analogue \eqref{aux_4} is the identity $$S(0) = \sum_{n\geq 0} G(n,0) = G(0,\frac12).$$
A detailed workflow of this example is available in the accompanying \textit{Mathematica} file.
\end{romexample}

\begin{remark}
For the classical-hypergeometric case, where an imaginary period is not available, Carlson's theorem is required in order to conclude that $S(k)$ is a constant.
\end{remark}

\begin{romexample}\label{ex_14}
The only other $1/\pi^4$-formula known, conjectured by Zhao and proved by the author \cite{au2025wilf}
	$$\sum _{n\geq 0} \left(-\frac{3^3}{2^8}\right)^n \frac{\left(\frac{1}{2}\right)_n^5 \left(\frac{1}{4}\right)_n \left(\frac{3}{4}\right)_n \left(\frac{1}{3}\right)_n \left(\frac{2}{3}\right)_n}{(1)_n^9}(4528 n^4+3180 n^3+972 n^2+147 n+9)=\frac{768}{\pi ^4}.$$
	It has a $q$-analogue
	$$\sum_{n\geq 0} \frac{(-1)^n q^{n^2}  P(n) (q;q^2)_{2 n} (q^2;q^2)_{3 n} (q;q^2)_n^{10}}{(q^{2 n}-q^n+1) (q^{2 n}+q^n+1) (q^{2 n+1}+1)^5 (q^2;q^2)_n^5 (q^2;q^2)_{2 n}^5}  = \frac{(q;q^2)_{\infty }^8}{(q^2;q^2)_{\infty }^8} $$
	where {\small $P(n) =1+q^{2 n} (1+5 q)+q^{1+6 n} (-5-16 q+10 q^2)+q^{4 n} (1-5 q+10 q^2)+q^{1+8 n} (-1-11 q-25 q^2+5 q^3)+q^{2+10 n} (19-14 q-15 q^2+q^3)+q^{12 n} (q^2+46 q^3-4 q^4-4 q^5)+q^{2+14 n} (-4-4 q+46 q^2+q^3)+q^{2+16 n} (1-15 q-14 q^2+19 q^3)-q^{3+18 n} (-5+25 q+11 q^2+q^3)+q^{4+20 n} (10-16 q-5 q^2)+q^{5+22 n} (10-5 q+q^2)+q^{6+24 n} (5+q)+q^{7+26 n}$}. 
To prove this, one simply applies the method in the example above to  $$\textsf{Dougall7F6}_{q^2}\left(\frac{1}{2}-n,\frac{1}{2}-n,\frac{1}{2}-n,\frac{1}{2}-n,\frac{1}{2}-n,k+n\right).$$
\end{romexample}

\begin{remark}
By putting accessory parameters into the WZ-seed, one can derive five-parameter extensions of both $q$-identities in Examples \ref{ex_13} and \ref{ex_14}. However, the full expressions are dauntingly long.
\end{remark}

\begin{romexample}\label{ex_15}
The $1/\pi^2$-formula $\sum_{n\geq 0} \left(\frac{1}{2^4}\right)^n \frac{(\frac{1}{2})_n^3 (\frac{1}{4})_n (\frac{3}{4})_n}{(1)_n^5} (120n^2 + 34n+3) = \frac{32}{\pi^2}$ in Example \ref{ex_9} has another $q$-analogue:
$$\sum_{n\geq 0} \frac{q^{2 n^2} P(n) (q;q^2)_n^2 (q;q^2)_{2 n}^2 (q^3;q^6)_n^2 (q^6;q^6)_{2 n}}{(1-q^{2 n}+q^{4 n}) (1+q^{1+2 n}) (1+q^{1+4 n}) (1+q^{2+4 n}) (q^2;q^2)_{2 n} (q^2;q^2)_{4 n} (q^6;q^6)_n^4} = \frac{(q;q^2)_{\infty }^2 (q^3;q^6)_{\infty }^2}{(q^2;q^2)_{\infty }^2 (q^6;q^6)_{\infty }^2},$$
\end{romexample}
here {\small $P(n) = 1+(q-1) q^{2 n}+q^{4 n} (1+q^2)+q^{3+6 n}+(-1+q) q^{1+8 n}+q^{2+10 n} (-2-q+q^2)-q^{3+12 n} (1+q)-q^{2+14 n} (1+q)+q^{2+16 n} (1-q-2 q^2)+(1-q) q^{4+18 n}+q^{3+20 n}+q^{22 n} (q^4+q^6)+(1-q) q^{5+24 n}+q^{6+26 n}$}. This, as well as a three-variable generalization\footnote{whose full expression can be found in the accompanying \textit{Mathematica} file} can be proved by starting from the WZ-seed
$$\textsf{Seed8}_q\left( \frac12+b-n,\frac12+b-2 n,c+n,k+n\right),$$
then proceed exactly as the previous two examples: $\sum_{n\geq 0} G(n,k)$ is a constant in $k$ whose value can be found by choosing a particular $k$ such that the sum terminates.

\begin{romexample}\label{ex_16}
The $1/\pi$-formula $$\sum_{n\geq 0} \left(\frac{1}{9}\right)^n \frac{\left(\frac{1}{2}\right)_n \left(\frac{1}{4}\right)_n \left(\frac{3}{4}\right)_n}{(1)_n^3} (8n+1) = \frac{2\sqrt{3}}{\pi}$$
has a $q$-analogue 
	$$\sum_{n\geq 0}\frac{q^{2 n^2} (1-q^{8n+1}) (q;q^2)_{2 n} (q;q^2)_n^2}{(q^2;q^2)_{2 n} (q^6;q^6)_n^2} =\frac{(q;q^2)_{\infty } (q^3;q^6)_{\infty }}{(q^2;q^2)_{\infty } (q^6;q^6)_{\infty }},$$
and a two-parameter extension:
\begin{equation}\label{aux_5}\sum_{n\geq0}\frac{q^{2bn+2cn+2n^2}(1-q^{b+4c+8n+1})(q^{b-2c+1},q^{b+4c+1};q^2)_n(q^{-b+2c+1};q^2)_{2n}}{(q^6;q^6)_n(q^{6c+6};q^6)_n(q^{2b+2c+2};q^2)_{2n}}=\frac{(q^{3b+3};q^6)_{\infty}(q^{b+4c+1};q^2)_{\infty}}{(q^{6c+6};q^6){}_\infty(q^{2b+2c+2};q^2)_{\infty}}.\end{equation}
To prove this, starting from $\textsf{Seed4}_{q^2}(1 - 2 n, \frac12+ b - 2 n, n + k)$, we are inspired to take
  \begin{multline*}F(n,k)=\frac{(1-q^{12k+6})(1-q^{6n-6k})(1-q^{2b+4n-2})(1-q^{2b+4n})}{(1-q^{6k})(1-q^{b-6k+2n-5})(1-q^{b-6k+2n-3})(1-q^{b-6k+2n-1})}\\\times\frac{q^{-3bk+2bn-9k+2n^2}(q^{1-b};q^2)_{2n}(q^6;q^6)_k(q^{b+1};q^2)_{3k}(q^{b-6k+1};q^2)_n(q^{b+6k+1};q^2)_n}{(q^{1-b};q^2)_{3k}(q^{6};q^6)_k(q^{2b+2};q^2)_{2n}(q^{6-6k},q^{6k+6};q^6)_n},\end{multline*}
then $G(n,k)$ equals $$ \frac{(1-q^{b+8n+1})q^{-3bk+2bn-9k+2n^2}(q^{1-b};q^2)_{2n}(q^6;q^6)_k(q^{b+1};q^2)_{3k}(q^{b-6k+1};q^2)_n(q^{b+6k+1};q^2)_n}{(1-q^{6k})(q^{1-b};q^2)_{3k}(q^{6};q^6)_k(q^{2b+3c+2};q^2)_{2n}(q^{6-6k},q^{6k+6};q^6)_n}.$$
Let $$F'(n,k) = F(n+k,k)+G(n+k,k+1), \quad G'(n,k) = G(n+k,k),$$
then $(F'(n,k), G'(n,k))$ is a WZ-pair since $(F(n,k),G(n,k))$ is. Also $F'(0,k) = 0$ and $\lim_{n\to\infty} F'(n,k) = 0$. Therefore the argument in the previous two examples applied to $(F',G')$ implies $S(k) := \sum_{n\geq 0} G'(n,k)$ is 1-periodic in $k$. That is\footnote{here we normalized $G'(n,k)$. The sum is still 1-periodic since we multiplied $G'(n,k)$ by $1$-periodic functions in $n,k$ only.}
$$\sum_{n\geq 0}\frac{q^{2bn+4kn+2n^2}(1-q^{b+8k+8n+1})(q^{b+1};q^2){}_{4k}(q^{b-4k+1};q^2)_n(q^{-b+4k+1};q^2)_{2n}(q^{b+8k+1};q^2)_n}{(q^6;q^6)_n(q^{6};q^6)_{2k}(q^{2b+2};q^2)_{2k}(q^{12k+6};q^6)_n(q^{2b+4k+2};q^2)_{2n}}$$
is 1-periodic in $k$. The above sum is moreover holomorphic in $k$; hence it is a constant. Letting $k=\frac{b-1}{4}$ gives its value. The formula \eqref{aux_5} is exactly this formula after renaming $k$ as $c$.
\end{romexample}

A one-parameter special case of the above was derived in \cite{guo2018ramanujan}.

\begin{romexample}\label{ex_17}
The remarkable $1/\pi^3$-formula
$$\sum_{n\geq 0} \left(\frac{1}{2^6}\right)^n \frac{\left(\frac{1}{2}\right)_n^7}{(1)_n^7} (6 n+1) (28 n^2+8 n+1) = \frac{32}{\pi^3}$$
has a $q$-analogue
\begin{equation}\label{aux_6}\sum_{n\geq 0}\frac{q^{2 n^2} (1-q^{6 n+1}) (1+q^{4 n}+q^{2 n+1}-6 q^{6 n+1}+q^{8 n+2}+q^{10 n+1}+q^{12 n+2}) (q^4;q^4)_{2 n} (q;q^2)_n^8}{(q^{4 n}+1) (q^{2 n+1}+1) (q^2;q^2)_{2 n}^2 (q^4;q^4)_n^6} = \frac{(q;q^2)_{\infty }^4}{(q^2;q^4)_{\infty }^2 (q^4;q^4)_{\infty }^6}.\end{equation}
\end{romexample}
To prove this, we start from the WZ-seed $$\textsf{Seed9}_{q^2}\left(1-2 n,-n+\frac{1}{2},-n+\frac{1}{2},k+n\right).$$
$$F(n,k) := \frac{(1-q^{8 k+4}) (1-q^{4 n})^4 (1-q^{4 n-2}) q^{2 n^2-4 k} (1-q^{4 n-4 k})^2 (q^4;q^4)_{2 n} (q;q^2)_n^4 (q^{1-4 k};q^2)_n^2 (q^{4 k+1};q^2)_n^2}{(1-q^{4 k})^2 (1-q^{8 n}) (1-q^{-4 k+2 n-3})^2 (1-q^{-4 k+2 n-1})^2 (q^2;q^2)_{2 n}^2 (q^4;q^4)_n^2 (q^{4-4 k};q^4)_n^2 (q^{4 k+4};q^4)_n^2},$$
then its WZ-mate is
$$G(n,k) = F(n,k) \frac{\splitfrac{q^{-4 k-2} (1-q^{6 n+1}) (q^{4 k+1}-q^{2 n})^2 (q^{4 k+3}-q^{2 n})^2 (q^{4 k+2 n+1}+q^{4 k+4 n}-2 q^{4 k+6 n+1}-2 q^{8 k+6 n+1}}{+q^{4 k+8 n+2}+q^{4 k+10 n+1}+q^{4 k+12 n+2}+q^{4 k}-2 q^{6 n+1})}}{(1-q^{8 k+4}) (1-q^{4 n})^3 (q^{4 n}-q^2) (q^{2 n+1}+1) (q^{4 k}-q^{4 n})^2}.$$
Similar to the above example, we let
$$F'(n,k) = F(n+k,k)+G(n+k,k+1), \quad G'(n,k) = G(n+k,k),$$
then $F'(0,k) = 0$ and $\lim_{n\to \infty} F'(n,k) = 0$, so $S(k) := \sum_{n\geq 0} G'(n,k)$ is 1-periodic in $k$ and holomorphic, so is a constant. Its value can then be found by letting $k=1/2$. The formula \eqref{aux_6} is exactly $\sum_{n\geq 0} G(n,0) = S(\frac{1}{2})$. \\
More generally, \eqref{aux_6} has a three-variable extension:
\begin{multline*}\sum_{n\geq 0} \frac{q^{b n+c n+2 d n+2 n^2} P(n) (q^{1+b-d},q^{1+c-d},q^{1-b+d},q^{1-c+d},q^{1+2 b+c+d},q^{1+b+2 c+d},q^{1+b+3 d},q^{1+c+3 d};q^2)_n (q^{4+2 b+2 c+4 d};q^4)_{2 n}}{\left(1-q^{2+b+c+2 d+4 n}\right) \left(1+q^{b+c+2 d+4 n}\right) (q^4,q^{4+2 b+2 c},q^{4+2 b+2 d},q^{4+2 b+2 c+4 d},q^{4+2 c+2 d},q^{4+4 d};q^4)_n (q^{2+b+c+2 d};q^2)_{2 n}^2}\\ = \frac{(q^{2+2 b+2 c},q^{2+4 c},q^{4+4 c};q^4)_{\infty } (q^{1+2 b+c+d},q^{1+b+2 c+d},q^{1+b+3 d},q^{1+c+3 d};q^2){}_{\infty }}{(q^{2+2 b+2 c},q^{2+4 c},q^{2+b+c+2 d},q^{2+b+c+2 d};q^2)_{\infty } (q^{4+2 b+2 d},q^{4+2 c+2 d},q^{4+4 d};q^4)_{\infty }},\end{multline*}
where {\small $P(n) = 1+q^{b+c+2 d+4 n} (1-q^2)-q^{1+d+6 n} (q^b+q^c) (q^{b+c}+q^{2 d}+q^{2 b+2 c+2 d}+q^{b+c+4 d})+q^{2+b+c+2 d+8 n} (1+q^{2 b+2 c}+q^{4 d}+q^{2 b+2 d}+q^{b+c+2 d}+q^{2 c+2 d}+q^{2 b+2 c+4 d})-q^{3+2 b+2 c+5 d+14 n} (q^b+q^c) (q^{b+c}+q^{2 d}+q^{2 b+2 c+2 d}+q^{b+c+4 d})+q^{2+2 b+2 c+4 d+12 n} (1+q^{2 b+2 c}+q^{4 d}+q^{2 b+2 d}+q^{b+c+2 d}+q^{2 c+2 d}+q^{2 b+2 c+4 d})-q^{2+4 b+4 c+8 d+16 n} (1-q^2)+q^{4+5 b+5 c+10 d+20 n}$}. \\
This is shown in exact parallel to the above by putting auxiliary parameters into the WZ-seed $$\textsf{Seed9}_{q^2}\left(1-2 n,b-n+\frac{1}{2},c-n+\frac{1}{2},k+n\right),$$
and then following the same recipe as in the above example. Renaming $k$ as $d$ produces the three-variable formula above. A detailed workflow of this example is available in the accompanying \textit{Mathematica} file.

\begin{romexample}\label{ex_18}
The $1/\pi$-formula $\sum_{n\geq 0} \left(\frac{1}{2^6}\right)^n \frac{(\frac{1}{2})_n^3}{(1)_n^3} (42n+5)= \frac{16}{\pi}$
has a more exotic $q$-analogue:
\begin{equation}\label{aux_7}\sum_{n\geq 0} \frac{q^{6 n^2} (q^{6 n+1}+q^{12 n+2}+1) P(n) (q^{12};q^{12})_{2 n} (q^3;q^6)_n^4 (q,q^5;q^6)_n^2}{(q^{12 n}+1) (q^{6 n+3}+1) (q^{6 n+1}+1)^2 (q^6;q^6)_{2 n}^2 (q^2,q^{10},q^{12};q^{12})_n^2} =  \frac{(q,q^5,q^7,q^{11};q^{12})_\infty (q^3,q^9;q^{12})_\infty^2}{(q^2,q^6,q^{10},q^{12};q^{12})_\infty^2},\end{equation}
where $P(n) = 1+q^{6 n} (q+q^3)+(q^{12 n}-q^{1+18 n}) (1+q^2+q^4)-q^{2+24 n} (1+q^2)-q^{5+30 n}.$ 
This follows from $$\textsf{Seed9}_{q^6}\left(\frac{5}{6}-2 n,-n+\frac{1}{2},-n+\frac{1}{6},k+n+\frac{1}{6}\right),$$
which gives
\begin{multline*}F(n,k) = \frac{(1-q^{24 k+14}) (1-q^{12 n})^2 (1-q^{12 n-6}) (1-q^{12 n-2}) q^{6 n^2-12 k} (1-q^{-12 k+12 n-2}) (1-q^{12 n-12 k})}{(1-q^{12 k}) (1-q^{24 n}) (1-q^{-12 k+6 n-11}) (1-q^{-12 k+6 n-9}) (1-q^{-12 k+6 n-5}) (1-q^{-12 k+6 n-3}) (1-q^{12 k+12 n+2})} \\ \times \frac{ (q^{12};q^{12})_{2 n} (q^{3-12 k},q^{1-12 k},q^{12 k+3},q^{12 k+5},q,q^3,q^3,q^5;q^6)_n}{ (q^6;q^6)_{2 n}^2 (q^2,q^{10};q^{12})_n (q^{10-12 k},q^{12-12 k},q^{12 k+12},q^{12 k+2};q^{12})_n}.\end{multline*}
The proof then proceeds \textit{mutatis mutandis} as in the previous two examples by constructing $F'(n,k), G'(n,k)$. \\
One can also derive a three-variable generalization\footnote{whose full expression can be found in the accompanying \textit{Mathematica} notebook} of \eqref{aux_7} by inserting accessory parameters into this WZ-seed.
\end{romexample}

\begin{romexample}\label{ex_19}
	The $1/\pi$-formula $\sum_{n\geq 0} \left(\frac{1}{4}\right)^n \frac{\left(\frac{1}{2}\right)_n^3}{(1)_n^3} (6n+1) = \frac{4}{\pi}$ has a $q$-analogue that is not modular, but mock modular as follows:
$$\sum_{n\geq 0}\frac{q^{2 n} (2-q^{2 n}-q^{4 n+1}) (q;q^2)_n^4}{(q^{2 n+1}+1) (q^2;q^2)_{2 n} (q^2;q^2)_n^2} = \frac{(q;q^2)_{\infty }^4}{(q^2;q^2)_{\infty }^3} \omega(q).$$
It has a three-variable extension:
\begin{multline}\label{aux_8}\sum_{n\geq 0}\frac{q^{2 n} (q^{b+c+2 d+6 n+2}+q^d (q^c-q^{2 n} (q^{b+c+1}+q^{b+1}+q^c))+q^b) (q^{b+1},q^{d+1},q^{b-c+1},q^{c+d+1};q^2)_n}{(1-q^{b+d+4 n+2}) (q^2;q^2)_n  (q^{-b+c+d+2};q^2)_n (q^{b+d+2};q^2)_{2 n}} \\  = \frac{q^{c-1} (q^{1+b},q^{1+b-c},q^{1+d},q^{1+c+d};q^2)_{\infty }}{(q^2,q^{2+b+d},q^{2-b+c+d};q^2)_{\infty }} \sum_{n\geq 1} \frac{q^{n(1+b-c)}}{(q^{1+b};q^2)_n},
\end{multline}
from which the above formula follows by setting $b=c=d=0$ and using $\omega(q) = \sum_{n\geq 1} \frac{q^{n-1}}{(q;q^2)_n}$. To prove this, we use the WZ-pair coming from $\textsf{Gauss2F1}_{q^2}(\frac{1}{2}-n,\frac{b}{2}+\frac{1}{2},\frac{c}{2}+1,-\frac{b}{2}+k-n+\frac{1}{2})$. That is,
\begin{multline*}F(n,k) = \frac{q^{-b k+c k+2 n} (1-q^{-b+c+2 n}) (1-q^{-2 k+2 n}) (1-q^{b-2 k+4 n}) }{(1-q^{2 k}) (1-q^{-b+2 k}) (1-q^{-1+b-2 k+2 n}) (1-q^{-1+b-c-2 k+2 n})} \\ \times \frac{ (q;q^2)_n (q^2;q^2)_k (q^{2-b};q^2)_k (q^{1+c};q^2)_n (q^{1+b-2 k};q^2)_n (q^{1+b-c-2 k};q^2)_n}{ (q^{1-b};q^2)_k (q^{1-b+c};q^2)_k (q^{2-b+c};q^2)_n (q^{2-2 k};q^2)_n (q^{2+b-2 k};q^2)_{2 n}},\end{multline*}
then its WZ-mate is \begin{multline*}G(n,k) =- \frac{q^{2-b k+c k+2 n} (q^{b+2 k}+q^{c+4 k}-q^{1+b+2 k+2 n}-q^{c+2 k+2 n}-q^{1+b+c+2 k+2 n}+q^{2+b+c+6 n})}{(1-q^{2 k}) (q^b-q^{2 k}) (q^{2 k}-q^{2+b+4 n})}\\ \times \frac{ (q;q^2)_n (q^2;q^2)_k (q^{2-b};q^2)_k (q^{1+c};q^2)_n (q^{1+b-2 k};q^2)_n (q^{1+b-c-2 k};q^2)_n}{ (q^{1-b};q^2)_k (q^{1-b+c};q^2)_k (q^{2-b+c};q^2)_n (q^{2-2 k};q^2)_n (q^{2+b-2 k};q^2)_{2 n}}. \end{multline*}
Analogous to the previous examples, let $$F'(n,k) = F(n+k,k)+G(n+k,k+1), \quad G'(n,k) = G(n+k,k),$$
then $F'(0,k) = 0$ and $\lim_{n\to \infty} F'(n,k) = 0$, so $\sum_{n\geq 0} G'(n,k)$ is 1-periodic in $k$ and holomorphic. As in the previous examples, this implies
$$S:=\sum_{n\geq 0} \frac{\splitfrac{q^{2 n} (q^b+q^{c+2 k}+q^{2+b+c+4 k+6 n}-q^{2 k+2 n} (q^{1+b}+q^c+q^{1+b+c}))}{ \times (q;q^2)_k (q^{1+c};q^2)_k (q^{1+2 k},q^{1+c+2 k},q^{1+b},q^{1+b-c};q^2)_n}}{(1-q^{2+b+2 k+4 n}) (q^{2+b};q^2)_k (q^{2-b+c};q^2)_k (q^{2+b+2 k};q^2)_{2 n} (q^2,q^{2-b+c+2 k};q^2)_n}$$
is independent of $k$. However, if one imitates the previous examples by choosing a particular $k$, then one encounters obstacles: every $k$ such that the above sum terminates is actually a singularity. For example, the sum terminates for $k=-1/2$ or $k=-(1+c)/2$. They are singularities for $(q;q^2)_k$ and $(q^{1+c};q^2)_k$ respectively. As no finite $k$ reduces the sum, we let $k\to \infty$, giving
$$S = \frac{q^b(q,q^{1+c};q^2)_\infty }{(q^{2+b},q^{2-b+c};q^2)_\infty }\sum_{n\geq 0} \frac{q^{2 n} (q^{1+b},q^{1+b-c};q^2)_n}{ (q^2;q^2)_n}.$$
One recognizes the sum on the RHS as the $q$-hypergeometric series  $_2\phi_1(q^{1+b},q^{1+b-c};0;q^2;q^2)$. Heine's transformation formula \cite[p.~13]{gasper2011basic}
$$_2\phi_1(a,b;c;q,z) = \frac{(b,az;q)_\infty}{(c,z;q)_\infty} {_2\phi_1}(c/b,z;az;q,b)$$
transforms it into $\sum_{n\geq 1} \frac{q^{n(1+b-c)}}{(q^{1+b};q^2)_n}$. Equating the two displayed expressions for $S$ and replacing $k$ by $d/2$ produces \eqref{aux_8}.
\end{romexample}

\begin{romexample}\label{ex_20}
	The $1/\pi^2$-formula $$\sum_{n\geq 0} \left(\frac{3^3}{4^3}\right)^n \frac{\left(\frac{1}{2}\right)_n^3 \left(\frac{1}{3}\right)_n \left(\frac{2}{3}\right)_n}{(1)_n^5} (74n^2 + 27n+3) = \frac{48}{\pi^2}$$ has a mock modular $q$-analogue
\begin{multline}\label{aux_9}\sum_{n\geq 0}\frac{q^{2 n} \left(3+3 q^{1+2 n}-2 (2-q) q^{1+4 n}-4 q^{1+6 n} (1+q)+2 (1-2 q) q^{1+8 n}+3 q^{2+10 n}+3 q^{3+12 n}\right) (q^2;q^2)_{3 n} (q;q^2)_n^6}{(q^{2 n}-q^n+1) (q^{2 n}+q^n+1) (q^{2 n+1}+1)^3 (q^2;q^2)_n^3 (q^2;q^2)_{2 n}^3} \\ = \frac{(q;q^2)_{\infty }^6}{(q^2;q^2)_{\infty }^5} \omega(q).
\end{multline}
To show this, we start from the WZ-seed $$\textsf{Balanced3F2}_{q^2}(\frac{1}{2}-n,\frac{1}{2}-n,n,1-n,k+2 n).$$
That is,
$$F(n,k) = \frac{q^{2 k+2 n} (1-q^{2 n})^2 (1-q^{4 n}) (q;q^2)_k^2 (q;q^2)_n^4 (q^{1+2 k};q^2)_n^2 (q^{2+2 k};q^2)_{3 n}}{(1-q^{2 k+6 n}) (q^2;q^2)_k^2 (q^2;q^2)_n^2 (q^2;q^2)_{2 n} (q^{2+2 k};q^2)_n (q^{2+2 k};q^2)_{2 n}^2},$$
its WZ-mate is
$$G(n,k) = \frac{q^{-2 k} P(n,k)}{(1-q^{2 n})^3 (1+q^{2 n}) (1+q^{1+2 n}) (1-q^{2+2 k+4 n})^2} F(n,k),$$
where {\small $P(n,k) = 2+q^{2 k}-q^{2+2 k+6 n} (-2+q-3 q^{2 k})-q^{5+6 k+18 n} (1+2 q^{2 k})+q^{3+6 k+14 n} (1+q^{1+2 k})+q^{4+6 k+16 n} (3+2 q+q^{1+2 k})-q^{2+4 k+12 n} (3 q-q^{2 k}+2 q^{1+2 k})-q^{2 n} (1+2 q^{2 k}+3 q^{1+2 k})-q^{4 n} (q+q^{2+2 k})-q^{1+4 k+10 n} (2+3 q+4 q^3+2 q^{1+2 k}+2 q^{3+2 k})+q^{1+2 k+8 n} (2+2 q^2+4 q^{2 k}+3 q^{2+2 k}+2 q^{3+2 k})$}. Since $F(0,k) = 0, \lim_{n\to\infty} F(n,k) = 0$, so (argue as in Example \ref{ex_13}) $\sum_{n\geq 0}G(n,k)$ is $1$-periodic in $k$, holomorphic and hence constant in $k$. Like the previous example, no finite $k$ reduces the sum, so we can only let $k\to \infty$, because $$\lim_{k\to\infty} G(n,k) = \frac{q^{2 n} (2-q^{2 n}-q^{1+4 n}) (q;q^2)_{\infty }^2 (q;q^2)_n^4}{(1+q^{1+2 n}) (q^2;q^2)_{\infty }^2 (q^2;q^2)_n^2 (q^2;q^2)_{2 n}},$$
we have $$\sum_{n\geq 0} G(n,0) = \frac{(q;q^2)_{\infty }^2}{(q^2;q^2)_{\infty }^2 } \sum_{n\geq 0}\frac{q^{2 n} (2-q^{2 n}-q^{1+4 n})  (q;q^2)_n^4}{(1+q^{1+2 n}) (q^2;q^2)_n^2 (q^2;q^2)_{2 n} }.$$
Here the RHS matches the LHS of equation \eqref{aux_8}, proving \eqref{aux_9}. Just like all other examples, there is a multiple-variable extension\footnote{which can be found in the accompanying \textit{Mathematica} notebook} of formula \ref{aux_9}.
\end{romexample}

\begin{romexample}\label{ex_21}
Finally we prove a mock modular $q$-analogue of $\sum_{n\geq 0} (\frac{1}{2^6})^n \frac{(\frac{1}{2})_n^3}{(1)_n^3} (42n+5)= \frac{16}{\pi}$:
\begin{equation}\label{aux_11}\sum_{n\geq 0} \frac{q^{2 n} (3+3 q^{1+2 n}-q^{4 n}+2 q^{2+4 n}-3 q^{1+6 n}-3 q^{2+8 n}-q^{3+10 n}) (q;q^2)_n^6}{(1+q^{1+2 n})^3 (q^2;q^2)_{2 n}^3} = -\frac{(q;q^2)_{\infty }^4}{q (q^2;q^2)_{\infty }^4} \sum_{n\in \mathbb{Z}} \frac{(-1)^n q^{n^2}}{1-q^{2n-1}}.\end{equation}
To prove this, start from the WZ-seed $\textsf{Gauss2F1}_{q^2}\left(\frac{1}{2}-n,\frac{1}{2}-n,1,k-n+\frac{1}{2}\right)$, let
$$F(n,k) = \frac{q^{2 n} (1-q^{4 n}) (1-q^{-2 k+4 n})^2 (q;q^2)_n^4 (q^2;q^2)_k^2 (q^{1-2 k};q^2)_n^2}{(1-q^{2 k})^2 (1-q^{-1-2 k+2 n})^2 (q;q^2)_k^2 (q^2;q^2)_{2 n} (q^{2-2 k};q^2)_{2 n}^2},$$
then
$$G(n,k) = F(n,k) \frac{\splitfrac{(q^{1+2 k}-q^{2 n})^2 (-q^{4 k}-2 q^{6 k}+3 q^{1+4 k+2 n}+q^{4 k+4 n}+q^{2+4 k+4 n}}{+q^{1+4 k+6 n}+q^{3+4 k+6 n}-2 q^{2+2 k+8 n}-2 q^{4+2 k+8 n}-2 q^{3+2 k+10 n}+q^{4+12 n}+q^{5+14 n})}}{(q^{2 k}-q^{4 n})^2 (1-q^{4 n}) (1+q^{1+2 n}) (q^{2 k}-q^{2+4 n})^2}.$$
Let
$$F''(n,k) = F(n+k,2k)+F(n+k,2k+1)+G(n+k,2k+2),\quad G''(n,k) = G(n+k,2k).$$
Then one easily checks that $(F'',G'')$ is a WZ-pair since $(F,G)$ is. Also $F''(0,k) = 0, \lim_{n\to \infty} F''(n,k) = 0$. Therefore $$\sum_{n\geq 0} G''(n,k) = \sum_{n\geq 0} G(n+k,2k)$$ is $1$-periodic in $k$. That is,
$$\sum_{n\geq 0} \frac{\splitfrac{q^{2 k^2-2 k+2 n} (1+2 q^{4 k}-3 q^{1+2 k+2 n}-q^{4 k+4 n}-q^{2+4 k+4 n}-q^{1+6 k+6 n}-q^{3+6 k+6 n}}{+2 q^{2+4 k+8 n}+2 q^{4+4 k+8 n}+2 q^{3+6 k+10 n}-q^{4+4 k+12 n}-q^{5+6 k+14 n}) (q;q^2)_k^2 (q^{1-2 k};q^2)_n^2 (q^{1+2 k};q^2)_n^4}}{(1-q^{1+2 n})^2 (1+q^{1+2 n})^2 (1+q^{1+2 k+2 n}) (q^2;q^2)_{2 k} (q^2;q^2)_{2 n}^2 (q^{2+4 k};q^2)_{2 n}},$$
is 1-periodic in $k$. It is also holomorphic, but unlike previous examples, it is not constant\footnote{recall the way we (Example \ref{ex_13}) deduce a summation to be constant in $k$ is to come up with an imaginary period $2\pi i /\log q$. In the current case, such an imaginary period does not exist due to the factor $q^{2k^2}$ in the numerator.} in $k$. Now we want to let $k\in \mathbb{N}$ and $k\to \infty$. This necessitates the splitting of $\sum_{n\geq 0}$ into $\sum_{n\geq k} + \sum_{0\leq n\leq k-1}$. These two parts give respective limits $$\sum_{n\geq 0} \frac{q^{2 n} (q;q^2)_{\infty }^4 (q;q^2)_n^2}{(q^2;q^2)_{\infty }^3}, \qquad \sum_{n\geq 0} \frac{q^{2 n+2 n^2} (q;q^2)_{\infty }^4}{(-1+q^{1+2 n})^2 (q^2;q^2)_{\infty }^3 (q;q^2)_n^2}.$$
Therefore $$\sum_{n\geq 0} G(n+k,2k) = \lim_{k\to  \infty} \sum_{n\geq 0} G(n+k,2k) = \frac{(q;q^2)_\infty^4}{(q;q^2)_\infty^3} \left(\sum_{n\geq 0} q^{2n}(q;q^2)_n^2 + \sum_{n\geq 0} \frac{q^{2n+2n^2}}{(q;q^2)_{n+1}^2}\right).$$
Using the following identity \cite[p.~67]{andrews2005ramanujan}:
$$\sum _{n\geq 1} q^{2 n} (q^{1-c};q^2)_n (q^{1-d};q^2)_n+\sum _{n\geq 0} \frac{q^{-2 n+c n+d n+2 n^2}}{(q^{1+c};q^2)_n (q^{1+d};q^2)_n} = -\frac{q^{-1+d} (q^{1-d};q^2)_{\infty }}{(q^2;q^2)_{\infty } (q^{1+c};q^2)_{\infty }} \sum _{n\in \mathbb{Z}} \frac{(-1)^n q^{c n+n^2}}{1-q^{-1+d+2 n}},$$
we see the two sums in parentheses combine to give $\sum_{n\in \mathbb{Z}}\frac{(-1)^n q^{n^2}}{1-q^{2n-1}}$, as claimed in \eqref{aux_11}.

By inserting accessory parameters in the WZ-seed above, the formula \eqref{aux_11} has the following three-parameter generalization:
\begin{multline*}\sum_{n\geq 0}\frac{q^{2 n} P(n) (q^{1-d},q^{1-2 b-d},q^{1-c-d},q^{1+2 b+d},q^{1+d},q^{1+2 b-c-d};q^2)_n}{(1-q^{2+4 n}) (q^{c+2 d}-q^{2+4 n}) (1-q^{2+2 b-c+4 n}) (q^2,q^{2+2 b-c},q^{2-c-2 d};q^2)_{2 n}} \\ = -q^{d-1}\frac{(q^{1-d},q^{1-c-d},q^{1+2 b-c-d},q^{1+d};q^2)_{\infty }}{(q^2,q^2,q^{2+2 b-c},q^{2-c-2 d};q^2)_{\infty }} \sum_{n\in \mathbb{Z}} \frac{(-1)^n q^{2 b n+d n+n^2}}{1-q^{d+2 n-1}} \end{multline*}
where {\small $P(n) = q^{2 b+2 d}+q^{c+2 d}+q^{2 b+c+4 d}-q^{1+d+2 n}-q^{1+2 b+d+2 n}-q^{1+2 b+3 d+2 n}-q^{1+4 b+3 d+2 n}-q^{1+c+3 d+2 n}-q^{1+2 b+c+3 d+2 n}+q^{2+2 d+4 n}-q^{2 b+2 d+4 n}+q^{2+4 b+2 d+4 n}+q^{2+2 b-c+8 n}+q^{4+2 b-c+8 n}+q^{2+2 b+2 d+8 n}+q^{4+2 b+2 d+8 n}+q^{2+4 b-c+2 d+8 n}+q^{4+4 b-c+2 d+8 n}-q^{5+2 b-c+d+10 n}-q^{5+4 b-c+d+10 n}-q^{4+4 b-2 c+12 n}-q^{4+2 b-c+12 n}-q^{4+4 b-c+2 d+12 n}+q^{6+4 b-2 c+16 n}$}.
\end{romexample}

The occurrence of mock modularity in the previous three examples seems puzzling, as it indicates a higher level of complexity than the previous modular examples. It might be worthwhile to look into more examples; here we give two more examples for interested readers to investigate further. Start with the WZ-seed
$$\textsf{Watson3F2}_{q^2}\left(\frac{1}{2}-n,n,\frac{1}{2}-n,\frac{1}{2}-n\right),$$
computing $G(n,0)$ indicates the following series
$$\sum_{n\geq 0} \frac{q^{2 n} P(n) (q;q^2)_n^4 (q^2;q^2)_{2 n} (q^2;q^4)_n (q^4;q^4)_{4 n}}{(1+q^{8 n}) (1+q^{1+2 n}) (1+q^{2+4 n})^3 (q^2;q^2)_{4 n} (q^4;q^4)_n (q^4;q^4)_{2 n}^3} $$
where {\small $P(n) =1+q^{1+2 n}+2 q^{3+6 n}+q^{3+10 n}-q^{4+16 n}-2 q^{4+20 n}-q^{6+24 n}-q^{7+26 n}+q^{4 n} (1+q^2)-q^{5+22 n} (1+q^2)+q^{2+8 n} (3+2 q^2)-q^{3+18 n} (2+3 q^2)+q^{2+12 n} (-2+q^2+q^4)+q^{1+14 n} (-1-q^2+2 q^4)$}, one expects that it also exhibits mock modularity. This would be a fifth $q$-analogue of $\sum_{n\geq 0} (\frac{1}{2^6})^n \frac{(\frac{1}{2})_n^3}{(1)_n^3} (42n+5)= \frac{16}{\pi}$. \par

The second example concerns a $q$-analogue of
$$\sum_{n\geq 0} \left(\frac{-3^6}{4^6}\right)^n \frac{\left(\frac12\right)_n \left(\frac13\right)_n \left(\frac23\right)_n \left(\frac16\right)_n \left(\frac56\right)_n}{(1)_n^5} (1930n^2+549n+45) = \frac{384}{\pi^2},$$
which itself was a long-standing conjecture \cite{guillera2016bilateral} until it was proved recently in \cite{au2025wilf}. A $q$-analogue could in principle be derived by considering
$$\text{Seed2}_{q^2}\left(\frac{1}{2}-n,n,1,2 n\right),$$
but analytical subtleties seem even more pronounced than in the preceding example.

\newpage

	\bibliographystyle{plain} % We choose the "plain" reference style
	\bibliography{../ref.bib} % Entries are in the refs.bib file
	
\end{document}